\documentclass[12pt, reqno]{amsart}
\usepackage{amssymb,latexsym,amsmath,amsfonts}
\usepackage{latexsym}
\usepackage[mathscr]{eucal}

8.7in \numberwithin{equation}{section}

\def ~{\hspace{1mm}}

\newtheorem{thm}{Theorem}[section]
\newtheorem{cor}[thm]{Corollary}
\newtheorem{lem}[thm]{Lemma}

\newtheorem{obs}[thm]{Observation}

\newtheorem{prop}[thm]{Proposition}
\newtheorem{defn}[thm]{Definition}
\newtheorem{rem}[thm]{Remark}
\newtheorem{note}[thm]{Note}
\numberwithin{equation}{section}

\def\textmatrix#1&#2\\#3&#4\\{\bigl({#1 \atop #3}\ {#2 \atop #4}\bigr)}
\def\dispmatrix#1&#2\\#3&#4\\{\left({#1 \atop #3}\ {#2 \atop #4}\right)}

\begin{document}
\title[Dilations of $\Gamma$-contractions]
{Dilations of $\Gamma$-contractions by solving \\operator equations}

\author[Bhattacharyya]{Tirthankar Bhattacharyya}
\address[Bhattacharyya]{Department of Mathematics, Indian Institute of Science, Banaglore 560 012}
\email{tirtha@member.ams.org}

\author[Pal]{Sourav Pal}
\address[Pal]{Department of Mathematics, Indian Institute of Science, Banaglore 560 012}
\email{sourav@math.iisc.ernet.in}

\author[ShyamRoy]{Subrata Shyam Roy}
\address[ShyamRoy]{Indian Institute of Science Education and Research, Pin code 741252, Nadia, West Bengal}
\email{ssroy@iiserkol.ac.in}

\thanks{The first author was supported in part by a grant from UGC via
DSA-SAP and UKIERI funded by the British Council and the second
author was supported by Research Fellowship of Council of Science
and Industrial Research, India and UKIERI.}

\begin{abstract}
For a contraction $P$ and a bounded commutant $S$ of $P$, we seek a
solution $X$ of the operator equation
$$S-S^*P=(I-P^*P)^{\frac{1}{2}}X(I-P^*P)^{\frac{1}{2}},$$ where
$X$ is a bounded operator on
$\overline{\mbox{Ran}}(I-P^*P)^{\frac{1}{2}}$ with numerical radius
of $X$ being not greater than $1$. A pair of bounded operators
$(S,P)$ which has the domain
$$\Gamma=\{(z_1+z_2,z_1z_2): |z_1|\leq1,|z_2|\leq1\}\subseteq \mathbb C^2$$ as a spectral set, is called a
$\Gamma$-contraction in the literature. We show the existence and
uniqueness of solution to the operator equation above for a
$\Gamma$-contraction $(S,P)$. This allows us to construct an
explicit $\Gamma$-isometric dilation of a $\Gamma$-contraction
$(S,P)$. We prove the other way too, i.e, for a commuting pair
$(S,P)$ with $\|P\|\leq1$ and the spectral radius of $S$ being not
greater than $2$, the existence of a solution to the above equation
implies that $(S,P)$ is a $\Gamma$-contraction. We show that for a
pure $\Gamma$-contraction $(S,P)$, there is a bounded operator $C$
with numerical radius not greater than $1$, such that $S=C+C^*P$.
Any $\Gamma$-isometry can be written in this form where $P$ now is
an isometry commuting with $C$ and $C^*$. Any $\Gamma$-unitary is of
this form as well with $P$ and $C$ being commuting unitaries.
Examples of $\Gamma$-contractions on reproducing kernel Hilbert
spaces and their $\Gamma$-isometric dilations are discussed.
\end{abstract}

\maketitle

\section{Motivation}

Subsets of $\mathbb{C}^n$ that are spectral sets or complete
spectral sets for a given commuting $n$-tuple of operators have been
studied for a long time, (see \cite{paulsen}) and the many
references cited there for the historical development.

Agler and Young in their seminal paper \cite{ay-jfa} introduced the
novel idea of studying all commuting pairs of bounded operators for
which a certain particular subset of $\mathbb{C}^2$ is a spectral
set. This subset is the symmetrized bidisc
 $$\Gamma=\{(z_{1}+z_{2},z_{1}z_{2}):|z_{1}|\leq1, |z_{2}|\leq1\}\subseteq \mathbb{C}^2 $$
 and the commuting pair of bounded operators $(S,P)$  defined on a Hilbert space $\mathcal{H}$ satisfies
$$\|f(S,P)\|\leq \sup_{(z_{1},z_{2})\in\Gamma} |f(z_{1},z_{2})|$$
where $f$ is a polynomial in two variables and the supremum is over
$\Gamma$. Thus, $\Gamma$ is a spectral set for $(S,P)$ or in other
words $(S,P)$ is a $\Gamma$- contraction. A $\Gamma$- contraction
$(S,P)$ is said to be a \em{pure} $\Gamma$-\em{contraction} if $P$
is a pure contraction, i.e, $P^{*n}\rightarrow 0$ as $n\rightarrow
\infty$. In other words $P\in C_{\cdot 0}$ following the terminology
of Sz-Nagy and Foias (see page-76 of \cite{nazy}). In their paper
(\cite{ay-jot}), Agler and Young described the motivation for
studying $\Gamma$-contractions. An understanding of this family of
operator pairs has led to the solutions of a special case of the
spectral Nevanlinna-Pick problem(\cite{ay-ieot}, \cite{ay-tams}),
which is one of the problems that arise in $H^{\infty}$ control
theory(\cite{BAF}). Also they play a pivotal role in the study of
complex geometry of the set $\Gamma$(see \cite{ay-blm}). In their
work Agler and Young did not assume separability of Hilbert spaces,
but in this note, all Hilbert spaces are over complex numbers and
are separable.

The remarkably smooth theory that they developed for
$\Gamma$-contractions parallels the highly successful theory of
dilation of a single contraction because they showed in
(\cite{ay-jot}), the existence of a $\Gamma$-isometric dilation for
any $\Gamma$-contraction. In this note we construct an explicit
$\Gamma$-isometric dilation of a $\Gamma$-contraction, i.e, given a
$\Gamma$-contraction $(S,P)$ on a Hilbert space $\mathcal H$, we
construct a space $\mathcal K$ containing $\mathcal H$ as a subspace
and a $\Gamma$-isometry $(T,V)$ on $\mathcal K$ such that
$T^*|_{\mathcal H}=S^*$ and $V^*|_{\mathcal H}=P^*$. In other words,
a $\Gamma$-contraction is the compression of a $\Gamma$-isometry to
a co-invariant subspace. What is remarkable here is that the space
$\mathcal K$ need not be any bigger than the minimal isometric
dilation space of the contraction $P$ and $V$ is in fact the minimal
isometric dilation of $P$. Moreover, $T$, in such a case, is
uniquely determined.

There are several ways to describe a $\Gamma$-contraction. We have
described a new way of characterizing $\Gamma$-contractions in
section 4. To do it, we define the fundamental equation of a pair of
bounded operators $(S,P)$ with $\|P\|\leq1$, to be the operator
equation
$$S-S^*P=\mathbf D_PX\mathbf D_P\,,\quad X\in\mathcal{B}(\mathcal
D_P).$$We show in the section on dilation, the existence and
uniqueness of solution to the fundamental equation for a
$\Gamma$-contraction $(S,P)$ and call the solution, the fundamental
operator for $(S,P)$. Uniqueness of minimal $\Gamma$-isometric
dilation (the minimality of a $\Gamma$-isometric dilation is defined
in section-2) of a $\Gamma$-contraction follows from the uniqueness
of the solution. This relates the theory of $\Gamma$-contractions
beautifully to solving operator equations. A one-parameter family of
examples of $\Gamma$-contractions has been obtained and is discussed
in section $3$. Their underlying spaces are reproducing kernel
Hilbert spaces. We give $\Gamma$-isometric dilations of those
$\Gamma$-contractions. Section $2$ describes the structure of
$\Gamma$-unitaries and $\Gamma$-isometries in complete detail with
some new characterizations of them.

We start by listing, without proof, some basic facts about the set
$\Gamma$ all of which can be found in \cite{ay-jot}. These will be
frequently used.

\begin{thm} Let $(s,p)\in \mathbb{C}^2$. The following are equivalent:
\begin{itemize}
\item[(i)] $(s,p)\in \Gamma$; \item[(ii)]
    $|s-\bar{s}p|+|p^2|\leq 1$ and $|s|\leq 2$; \item[(iii)]
    $2|s-\bar{s}p|+|s^2-4p|+|s^2|\leq 4$;
    \item[(iv)]$\rho(\alpha s,{\alpha}^2p)\geq 0$ , for all $
    \alpha\in\mathbb{D}$, where $\mathbb{D}$ is the unit open
    disc in $\mathbb{C}$; \item[(v)]$|p|\leq 1$ and  there
    exists $\beta\in\mathbb{C}$ such that $|\beta|\leq1$ and
    $s = \beta +\bar{\beta} p$;
    \item[(vi)]$|s|\leq2$ and $|(2\alpha p-s)(2-\alpha
        s)^{-1}|\leq1$\; \mbox{for all }\;
        $\alpha\in\mathbb{D}$;
    \item[(vii)]$1-\bar{\alpha}s+ \bar{\alpha}^2p\neq0$ and
        $|(p-\alpha
        s+{\alpha}^2)(1-\bar{\alpha}s+\bar{\alpha}^2p)^{-1}|\leq1$\;
        \mbox{for all }\; $\alpha\in\mathbb{D}$ .
\end{itemize}
\end{thm}

\begin{defn}
The $distinguished \;boundary$ of the set $\Gamma,$ denoted by
$b\Gamma$ is defined to be the set
$$b\Gamma=\;\{(z_{1}+z_{2}\;,z_{1}z_{2}):\; |z_1|=|z_{2}|=1\}.$$
\end{defn}

This is the \u{S}ilov boundary of the algebra of functions
continuous on $\Gamma$ and analytic in the interior of $\Gamma$.
\begin{thm} \label{1.3}
Let $(s,p)\in\mathbb{C}^2.$ Then the following are equivalent:
\begin{enumerate}
\item $(s,p)\in b\Gamma ;$ \item
    $(s,p)=(2xe^{i{\frac{\theta}{2}}},e^{i{ \theta}})$ for
    some $\theta\in\mathbb{R},$ and $x\in [-1,1].$
    \item $|p|=1,\; s=\overline{s}p \; and \;|s|\leq2$.
    \item $|p|=1,\;s = \beta  + \bar{\beta} p$ for some
        $\beta\in\mathbb{C}$ of modulus $1$.
\end{enumerate}

\end{thm}
We give a proof of $(1)\Leftrightarrow (4)$ because we could not
locate it in literature.
\begin{proof}
Let $|p|=1$ and $s=\beta +\bar{\beta}p$ for some
$\beta\in\mathbb{C}$ of modulus $1$. Taking $z_1=\beta$ and
$z_2=\bar{\beta}p$ we see that $s=z_1+z_2$ and $p=z_1z_2$ where
clearly $|z_1|=|z_2|=1$. Hence $(s,p)\in b\Gamma$.\\
Conversely, let $(s,p)\in b\Gamma$. Then $s=z_1+z_2$ and $p=z_1z_2$
for some $z_1,z_2$ of modulus $1$. Clearly $|p|=1$ and
$z_2=\bar{z_1}p$. Thus we have $s=z_1+
\bar{z_1}p=\beta+\bar{\beta}p$, where $\beta=z_1$.
\end{proof}

\begin{lem} \label{polcon}
$\Gamma$ is polynomially convex but not convex.
\end{lem}
There are more results about $\Gamma$ that we are not going into
because those are not relevant here. Symmetrized polydisc has been
studied in detail. The interested reader is referred to
\cite{ay-ems}, \cite{ay-blm}, \cite{ay-jga}, \cite{edi-zwo}.

\section{Structure theorems for $\Gamma-$isometries and $\Gamma-$unitaries}

Ever since Sz.-Nagy found the minimal unitary dilation for a
contraction on a Hilbert space, it became clear how powerful a tool
it is for studying an arbitrary contraction. An operator $T$ is a
contraction if and only if $\|p(T)\|\leq \|p\|_\infty$ for all
polynomials $p$ by von Neumann's inequality. This property can be
isolated and a compact subset $X$ of $\mathbb{C}$ is called a
spectral set for an operator $T$ if \begin{eqnarray}\label{a0} \|
f(T)\|\leq \sup_{z\in X} \|f(z)\|
\end{eqnarray}
for all rational functions $f(z)$ with poles off $X$ (we bring in
rational functions instead of just polynomials because the domain
$X$ is assumed to be just compact and not necessarily simply
connected, unlike $\mathbb D$). If (\ref{a0}) holds for all matrix
valued rational functions $f$, then $X$ is called a complete
spectral set for $T$. Moreover, $T$ is said to have a normal
$\partial X$-dilation if there is a Hilbert space $\mathcal{K}$
containing $\mathcal{H}$ as a subspace and a normal operator $N$ on
$\mathcal{K}$ with $\sigma(N)\subseteq
\partial X$ such that $$f(T)=P_{\mathcal H}f(N)|_{\mathcal H},$$
for all rational functions $f$ with poles off $X$. It is a
remarkable consequence of Arveson's extension theorem that $X$ is a
complete spectral set for $T$ if and only if $T$ has a normal
$\partial X$-dilation. Rephrased in this language, Sz.-Nagy dilation
theorem says that if $\mathbb{D}$ is a spectral set for $T$ then $T$
has a normal $\partial \mathbb{D}$-dilation. For $T$ to have a
normal $\partial X$-dilation it is necessary that  $X$ be a spectral
set for $T$. Sufficiency has been investigated for many domains in
$\mathbb{C}$ and several interesting results are known including
failure of such a dilation in multiply connected domains \cite{DM}.
If $X\subseteq \mathbb{C}^2$, then the questions are much more
subtle. If $(T_1,T_2)$ is a commuting pair of operators for which
$\mathbb{D}^2$ is a spectral set, then $(T_1,T_2)$ has a
simultaneous commuting unitary dilation by Ando's theorem. Taking
cue from such classically beautiful concepts, Agler and Young
introduced the following definitions.

\begin{defn} A commuting pair $(S,P)$ is called a $\Gamma$-{\em unitary}
if $S$ and $P$ are normal operators and the joint spectrum
$\sigma(S,P)$ of $(S,P)$ is contained in the distinguished boundary
of $\Gamma$. \end{defn}

\begin{defn}
A commuting pair $(S,P)$ is called a $\Gamma$-isometry if there
exist a Hilbert space $\mathcal{N}$ containing $\mathcal{H}$ and a
$\Gamma$-unitary $(\tilde{S},\tilde{P})$ on $\mathcal{N}$ such that
$\mathcal{H}$ is left invariant by both $\tilde{S}$ and $\tilde{P}$,
and
$$S = \tilde{S}|_{\mathcal{H}} \mbox{ and } P = \tilde{P}|_{\mathcal{H}}.$$

In other words, $(\tilde{S},\tilde{P})$ is a $\Gamma$-unitary
extension of $(S,P)$. A commuting pair $(S,P)$ is a $\Gamma$-{\em
co-isometry} if $(S^*,P^*)$ is a $\Gamma$-isometry. Moreover, a
$\Gamma$-isometry $(S,P)$ is said to be a pure $\Gamma$-isometry if
$P$ is a pure isometry, i.e, there is no non trivial subspace of
$\mathcal H$ on which $P$ acts as a unitary operator.\end{defn}

Here and henceforth, when we say joint spectrum, we shall mean the
Taylor joint spectrum unless otherwise mentioned. Let
\begin{align*} \rho( S,P)&=\;
2(I-P^*P)-(S-S^*P)-(S^*-P^*S)\\
&=\; \frac{1}{2}\{(2-S)^*(2-S)\;-(2P-S)^*(2P-S)\}.
\end{align*} The following result was proved in \cite{ay-jfa}.
There, in fact, it was proved that positivity $ \rho( S,P) $ is a
necessary and sufficient condition for $(S,P)$ to be a
$\Gamma$-contraction. A straightforward proof of one direction is
given below using joint spectral theory. Stinespring dilation is
avoided for proving this because this result will be used for
constructing explicit dilations.

\begin{prop} \label{prop1}
Let $(S,P)$ be a $\Gamma$-contraction. Then $\;\rho(\alpha
S,{\alpha}^2 P)\;\geq0,$ for all $\alpha\in\overline{\mathbb{D}}$.
\end{prop}

\begin{proof}

Let $\sigma(S,P)$ denote the Taylor joint spectrum of $(S,P)$. By
Lemma 6.11 of Chapter-III of \cite{vasilescu},
$$\sigma(S,P)\subset\;\sigma_{\mathcal{U}}(S,P),$$ where
$\mathcal{U}$ is the \textit{Banach subalgebra} of
$\mathcal{B}(\mathcal{H}),$ generated by $S,P$ and $I$ and
$\sigma_{\mathcal{U}}(S,P)$ is the joint spectrum of $(S,P)$
relative to this commutative Banach algebra.

It is straightforward from the definition of $\Gamma$ contraction
that
$$\; \sigma_{\mathcal{U}}(S,P)\subseteq \Gamma,$$
and hence we have $\sigma(S,P)\subseteq\Gamma.$

Let $f$ be a holomorphic function in a neighbourhood of $\Gamma.$
Since $\Gamma$ is polynomially convex, by Oka-Weil theorem (Theorem
5.1 of \cite{oka-weil}) there exists a sequence of polynomials
$\{p_{n}\}$ that converges uniformly to $f$ on $\Gamma$. So by
Theorem 9.9 of Chapter-III of \cite{vasilescu} we have
$$p_{n}(S,P)\rightarrow\;f(S,P)$$ which by virtue of $(S,P)$ being a $\Gamma$-contraction implies that
$$ \|f(S,P)\| = \lim_{n\rightarrow\infty} \|p_{n}(S,P)\|\leq \lim_{n\rightarrow\infty} \|p_{n}\|_{\Gamma}
=\|f\|.$$ Using the function $f(s,p)=({2{\alpha}^2 p-\alpha
s})/({2-\alpha s})$ which is holomorphic in a neighbourhood of
$\Gamma$ for $\alpha \in \mathbb{D},$ we get
$$\|\,(2{\alpha}^2 P - \alpha S)(2- \alpha S)^{-1}\,\|\leq\;\|f\|_{\Gamma}\,\leq1.$$
Thus ${(2-\alpha S)^*}^{-1}(2{\alpha}^2P-\alpha
S)^*(2{\alpha}^2P-\alpha S)(2-\alpha S)^{-1}\leq I.$\\This happens
if and only if $(2-\alpha S)^*(2-\alpha S)\geq (2{\alpha}^2P-\alpha
S)^*(2{\alpha}^2P-\alpha S).$ By definition of $\rho (S,P)$, the
last inequality is the same as $\rho (\alpha S,{\alpha}^2P)\geq 0.$

By continuity, $\rho(\alpha S, {\alpha}^2 P)\geq 0$ for all $ \alpha
\in \overline{\mathbb{D}}$.
\end{proof}

It is clear from the definition that if $(S_1, P_1)$ and $(S_2,
P_2)$ are $\Gamma$-unitaries, then so is the direct sum $(S,P) =
(S_1 \oplus S_2 , P_1 \oplus P_2)$. Indeed, the joint spectrum of
$(S,P)$ is the union of the joint spectrum of $(S_1, P_1)$ and the
joint spectrum of $(S_2, P_2)$ (see \cite{curto}). We begin with an
elementary lemma whose proof we skip because it is routine.

\begin{lem} \label{basiclemma} Let $X$ be a bounded operator on a Hilbert space
$\mathcal H$. If Re $\beta X \le 0$ for all complex numbers $\beta$
of
    modulus $1$, then $X = 0$. \end{lem}

Parts of the following theorem, which gives new characterizations of
$\Gamma$-unitaries were obtained by Agler and Young in
(\cite{ay-jot}). Parts (3), (4) and (5) are new.

\begin{thm} \label{G-unitary}
Let $(S,P)$ be a pair of commuting operators defined on a Hilbert
space $\mathcal{H}.$ Then the following are equivalent:
\begin{enumerate}

\item $(S,P)$ is a $\Gamma$-unitary\; ; \item there exist
    commuting unitary operators $U_{1}$ and $
    U_{2}$ on $\mathcal{H}$ such that
$$S= U_{1}+U_{2},\quad P= U_{1}U_{2}\; ;$$
\item $P$ is unitary,\;$S=S^*P,\;$\;and $r(S)\leq2,$ \; where
    $r(S)$ is the spectral radius of $S$.
    \item $(S,P)$ is a $\Gamma$-contraction and $P$ is a
        unitary.
    \item $P$ is a unitary and $S = U + U^* P$ for some
        unitary $U$ commuting with $P$.
\end{enumerate}
\end{thm}
\begin{rem}
We draw attention to the similarity between part$(5)$ of this
theorem and part$(4)$ of Theorem 1.3.
\end{rem}

\begin{proof}
$(1)\Rightarrow (2)$ This proof is the same as the one given by
Agler and Young in \cite{ay-jot}. We include it for the sake of
completeness. Let $(S,P)$ be a $\Gamma$-unitary. By the spectral
theorem for commuting normal operators there exists a
spectral measure say $M(.)$ on $\sigma(S,P)$ such that \\
$$ S=\; \int_{\sigma(S,P)}\; p_{1}(z)M(dz), \quad P=\; \int_{\sigma(S,P)}\;
p_{2}(z)M(dz),$$where $p_{1},\;p_{2}$ are the co-ordinate functions
on $\mathbb{C}^2.$ Now choose a measurable right inverse $\beta$ of
the restriction of the function $\pi$ to $\mathbb{T}^2$ so that
$\beta$ maps the distinguished boundary $b\Gamma$ of $\Gamma$ to
$\mathbb{T}^2.$ Let $\beta=\;(\beta_{1},\beta_{2})$ and $$U_{j}=\;
\int_{\sigma(S,P)}\; \beta_{j}(z)M(dz),\quad j=1,2.$$ Then $U_{1},
U_{2}$ are commuting unitary operators on $\mathcal{H}$ and
$$U_{1}+U_{2}=\;
\int_{\sigma(S,P)}\;(\beta_{1}+\beta_{2})(z)M(dz)\;=\;
\int_{\sigma(S,P)}\; p_{1}(z)M(dz)\;=\;S.$$ Similarly
$U_{1}U_{2}=\; P.$ Thus $(1)\Rightarrow(2).$\\
$(2)\Rightarrow(3)$ is clear.\\
$(3)\Rightarrow(1)$ We have $P^*P=PP^*=I$ and $S=S^*P$. Therefore,
$S^*=P^*S$ and as a consequence $$SS^*=(S^*P)(P^*S)=S^*S$$ as $P$ is
unitary. So $(S,P)$ is a commuting pair of normal operators. So
we have $r(S)=\|S\|$.\\
 Let $C^*(S,P)$ be the commutative
$C^*$-algebra generated by them. By general theory of joint spectrum
(see p-27, Proposition 1.2 of \cite{curto}),
$$\sigma(S,P) = \{(\varphi(S)\;,\varphi(P)): \varphi \in\mathcal{M}\},$$ where
$\mathcal{M}$ is the maximal ideal space of $C^*(S,P).$ Let
$(s,p)=(\psi(S),\psi(P))\in\sigma(S,P),$ where $\psi\in\mathcal{M}.$
Then
$$|p|^2 = \overline{p}p = \overline{\psi(p)}\psi(p) = \psi(P^*)\psi(P) = \psi(P^*P) = \psi(I)=1$$
and
$$\overline{s}p = \overline{\psi(S)}\psi(P) = \psi(S^*P) = \psi(S) = s.$$
Also $|s| =|\psi(S)| \le \|S\|=r(S) \leq 2.$ Therefore, by Theorem
\ref{1.3}, $(s,p)\in b\Gamma$ i.e, $\sigma(S,P)\subseteq b\Gamma$.
So $(S,P)$ is a $\Gamma$-unitary. Hence (3) $\Rightarrow$ (1).

Thus (1), (2) and (3) are equivalent.

The implication (1) $\Rightarrow$ (4) is trivial.

(4) $\Rightarrow$ (3) depends on the fact that if $(S,P)$ is a
$\Gamma$-contraction, then
$$\rho(\alpha
S,{\alpha}^2P)\geq0\;,\quad \mbox{ for all }
\alpha\in\overline{\mathbb{D}}.$$ Therefore, for $\beta \in
\mathbb{T}$, we have $ \rho(\beta S, {\beta}^2P) = 2(I-P^*P) -
{\beta}(S-S^*P) - {\overline{\beta}}(S^*-P^*S)\geq 0$. Using the
fact that $P^*P=I$, we get that  Re $\beta(S-S^*P)\leq 0$. By
invoking Lemma \ref{basiclemma} now, we get that $S-S^*P=0$.  Also
since $(S,P)$ is a $\Gamma$-contraction, $r(S)\leq \|S\|\leq2.$
Hence done.

(2) $\Rightarrow$ (5) follows as $S = U_1 + U_2 = U_1 + U_1^* P$ and
$U_1 P = U_1 U_1 U_2 = U_1 U_2 U_1 = P U_1$.

(5) $\Rightarrow$ (2) follows by taking $U_1 = U$ and $U_2 = U^* P$.

\end{proof}

\begin{cor}
The pair $(S,I)$ can be a $\Gamma$-contraction only by being a
$\Gamma$-unitary. It is so if and only if $S$ is a self-adjoint
operator of spectral radius not bigger than $2$.\end{cor}

During the course of the proof, we used something which we segregate
as a separate result because it will be used later too.

\begin{obs} If $P$ is a unitary, $S$ commutes with $P$ and
$S=S^*P$, then $S$ is normal. \label{unitary} \end{obs}

The structure of $\Gamma$-isometries can be deciphered using
numerical radius. We recall the definitions of $numerical$ $range$
and $numerical$ $radius$ and discuss some of their properties which
will be useful. The $numerical$ $range$ of an operator $T$ on a
Hilbert space $\mathcal{H}$ is defined to be
$$\Omega(T) =  \{\langle Tx,x \rangle\; : \; \|x\|_{\mathcal{H}}\leq 1\}.$$ The $numerical$ $radius$ of $T$ is
defined as
$$\omega(T) = \sup \{|\langle Tx,x \rangle|\; : \;
\|x\|_{\mathcal{H}}\leq 1\}.$$  It is well known that $r(T)\leq
\omega(T)\leq \|T\|$ for a bounded operator $T$. An elementary fact
will be used more than once, and hence we state it as a lemma
followed by a remarkable result due to Ando.

\begin{lem} \label{basicnrlemma} The numerical radius of an operator $X$ is not greater than one
if and only if  Re $\beta X \le I$ for all complex numbers $\beta$
of modulus $1$.
\end{lem}

\begin{proof}
It is obvious that $\omega(X)\leq1$ implies that $\mbox{Re }\beta
X\leq I$ for all $\beta\in\mathbb T$. We prove the other way. By
hypothesis, $\langle$Re $ \beta X h,h\rangle \leq1$ for all
$h\in\mathcal{H}$ with $\|h\|\leq1$ and for all $\beta \in
\mathbb{T}$. Note that $\langle$Re $ \beta X h,h\rangle =$ Re $
\beta \langle  X h,h\rangle$. Write $\langle Xh,h \rangle = e^{i
\varphi_h}|\langle Xh,h\rangle|$ for some $\varphi_h \in\mathbb{R},$
and then choose $\beta = e^{- i \varphi_h}$. Then we get  $|\langle
Xh,h \rangle|\leq1$ and this holds for each $h\in\mathcal{H}$ with
$\|h\|\leq1.$ Hence done. \end{proof}

\begin{thm} \label{ando}
\textbf{(Ando):} The numerical radius of an operator $X$ is not
greater than one if and only if there is a contraction $C$ such that
$$X=2(I-C^*C)^{1/2}C.$$
\end{thm}
For details of the proof, see Theorem 2 of \cite{ando}.

\begin{defn}
A bounded operator $X$ is said to be hyponormal if $X^*X \geq XX^*.$
\end{defn}
\begin{prop}
\textbf{(Stampfli):} \label{stampfli}If $X$ is hyponormal, then
$\|X^n\|=\|X\|^n$ and so $\|X\|=r(X)$.
\end{prop}
For details of the proof see Proposition 4.6 of \cite{conway}.
\begin{lem}\label{123}
Let $U\;,\;V$ be a unitary and a pure isometry on Hilbert Spaces
$\mathcal{H}_1 \;,\; \mathcal{H}_2$ respectively, and let
$X\;:\;\mathcal{H}_1 \rightarrow \mathcal{H}_2$ be such that
$XU=VX.$ Then $X=0.$
 \end{lem}
\begin{proof}
We have , for any positive integer $n,\; XU^n\;=\;V^nX$ by
iteration. Therefore, ${U^*}^nX^*\;=\;{X^*}{V^*}^n.$ Thus $X^*$
vanishes on $Ker{V^*}^n,$ and since $\displaystyle
\bigcup_n{KerV^*}^n$ is dense in $\mathcal{H}_2$ we have $X^*=0\;
i.e,\; X=0.$
\end{proof}

\begin{thm} \label{Gamma-isometry}
Let $S,P$ be commuting operators on a Hilbert space $\mathcal{H}.$
The following statements are all equivalent:\begin{enumerate}
\item $(S,P)$ is a $\Gamma$-isometry, \item $(S,P)$ is a
$\Gamma$-contraction and $P$ is isometry, \item
    $P$ is an isometry\;,\;$S=S^*P$ and $r(S)\leq2,$ \item
    $r(S)\leq2$ and $\rho(\beta S,{\beta^2P})=0$ for all
    $\beta\in\mathbb{T}$.

Moreover if the spectral radius $r(S)$ of $S$ is less than $2$ then
{\em (1)},{\em (2)},{\em (3)} and {\em (4)} are equivalent to: \item
$(2{\beta}P-S)(2-\beta S)^{-1}$ is an isometry, for all
$\beta\in\mathbb{T}.$
\end{enumerate}
\end{thm}

\begin{proof}
(1)$\Rightarrow$(2) is obvious.

(2)$\Rightarrow$(3) The fact that $(S,P)$ is a $\Gamma$-contraction
implies that $\| S \| \le 2$, whence $r(S) \le 2$. It also implies
that $ \rho (\alpha S , \alpha^2 P) \ge 0$ for all $\alpha$ in the
closed disk, in particular on the circle. In view of $P$ being an
isometry, this means that
$$ \mbox{ Re } \beta (S - S^* P) \le 0$$
for all $\beta $ of modulus $1$. By using Lemma \ref{basiclemma}, we
get that $S = S^* P$.

(3)$\Rightarrow$(4) This is obvious.

(4)$\Rightarrow$(1)  We have
$$\rho(\beta S,{\beta}^2
P) = 2(I-P^*P)- \beta(S-S^*P)-{\overline{\beta}}(S^*-P^*S)=0 \mbox{
for all } \beta\in\mathbb{T}.$$ Putting $\beta=1$ and $\beta=-1$, we
get $P^*P=I$ from which it follows by the same argument as above
that $S=S^*P.$ We shall now show that $(S,P)$ is a $\Gamma$-isometry
by exhibiting a $\Gamma$-unitary extension.

Wold decomposition of the isometry $P$ breaks the whole space
$\mathcal H$ into the direct sum of two reducing subspaces
${\mathcal H}_1$ and ${\mathcal H}_2$ so that $P$ has the form
$$P = \begin{pmatrix}P_1&0\\0&P_2\end{pmatrix}\mbox{ on }  {\mathcal{H}_1}\oplus {\mathcal{H}_2}=
\mathcal{H},$$ where $P_1$ is a unitary and $P_2$ is a pure isometry
(a shift of some multiplicity). With respect to this decomposition
of $\mathcal H$, we write
$$S\;=\;
\begin{pmatrix}S_{11}&S_{12}\\S_{21}&S_{22}\end{pmatrix}.$$ By commutativity of $S$ and $P$ and applying
Lemma \ref{123}, we see that $S$ takes the form
$$\begin{pmatrix}S_1&0\\0&S_2\end{pmatrix}\quad \mbox{on}\; \mathcal{H}_1
\oplus \mathcal{H}_2\quad \mbox{and} \quad S_1P_1=P_1S_1\;,\;
S_2P_2=P_2S_2.$$ Also by $S=S^*P$ and $P^*P=I$ we get
$S_i=S_i^*P_i$ and $P_i^*P_i=I$ respectively for $i=1,2.$\\
The pair $(S_1,P_1)$ is $\Gamma$-unitary by part (4) of Theorem
\ref{G-unitary} because $P_1$ is unitary and restriction of a
$\Gamma$-contraction to an invariant subspace is a
$\Gamma$-contraction.\\
 The pair $(S_2,P_2)$ is a $\Gamma$-contraction since $(S,P)$ is so. Since $P_2$ is a pure isometry it
can be identified with the multiplication operator $M_z^E$ on
$H^2(E)$ for some Hilbert space $E.$ Again since $S_2$ commutes with
$P_2(\equiv M_z^E)$, it can be identified with the multiplication
operator $M_{\varphi}^E$ for some $\varphi\in
{H^{\infty}(\mathcal{B}(E))}.$\\
Also because $P_2$ is isometry, $I-P_2P_2^*\geq 0$ and we have

$$S_2^*(I-P_2P_2^*)S_2\geq0 \Rightarrow S_2^*S_2 \geq (S_2^*P_2)(P_2^*S_2)=S_2S_2^*, \mbox{ since }S_2=S_2^*P_2.$$ Thus $S_2$ is
hyponormal and by Stampfli's result (Theorem \ref{stampfli}),
$r(S_2)=\|S_2\|$ and hence
$\|\varphi\|=\|M_{\varphi}^E\|=\|S_2\|\leq2.$ Since
$S_2={S_2}^*P_2$, or equivalently
$M_{\varphi}^E={M_{\varphi}^E}^*M_z^E$, we have
$$\varphi(z)={\varphi}^*(z)z \mbox{ for all } z\in\mathbb{T}.$$
Consider on $L^2(E)$, the multiplication operators $U_{\varphi}^E$
and $U_z^E$, multiplication by $\varphi(z)$ and $z$ respectively.
Obviously $U_z^E$ is a unitary operator on $L^2(E).$ Since
$\varphi(z)={\varphi}^*(z)z$ we have
$U_{\varphi}^E={U_{\varphi}^E}^*U_z^E$, i.e,
${U_{\varphi}^E}^*={U_z^E}^*U_{\varphi}^E$ and hence
$$U_{\varphi}^E{U_{\varphi}^E}^*=({U_{\varphi}^E}^*U_z^E)({U_z^E}^*{U_{\varphi}}^E)={U_{\varphi}^E}^*U_{\varphi}^E$$
and thus $U_{\varphi}^E$ is normal. So we have a pair of commuting
normal operators $(U_{\varphi}^E,U_z^E)$ on $L^2(E))$ such that
$r(U_{\varphi}^E)=\|U_{\varphi}^E\|=\|{\varphi}\|\leq2\;,\;U_{\varphi}^E={U_{\varphi}^E}^*U_z^E$
and $U_z^E$ is unitary. Therefore by part (3) of Theorem
\ref{G-unitary}, $(U_{\varphi}^E,U_z^E)$ is a $\Gamma$-unitary. The
restriction to $H^2(E)$ of this $\Gamma$-unitary is
$(M_{\varphi}^E,M_z^E).$ In other words $(U_{\varphi}^E,U_z^E)$ is a
$\Gamma$-unitary extension of $(M_{\varphi}^E,M_z^E)$.

Taking $\tilde{S}= S_1 \oplus U_{\varphi}^E$  and $\tilde{P}=P_1
\oplus U_z^E$ on ${\mathcal H}_1 \oplus L^2(E)$, we see that
$(\tilde{S},\tilde{P})$ is a $\Gamma$-unitary extension of $(S,P)$.
Hence $(S,P)$ is a $\Gamma$-isometry.

Thus (1) through (4) are equivalent.

(4) $\Leftrightarrow$ (5) By hypothesis,
\begin{align*} & \rho(\beta S,{\beta}^2P)=\;\frac{1}{2}\{(2-\beta
S)^*(2-\beta S)\;-\;(2{\beta}^2P-\beta S)^*(2{\beta}^2P-\beta
S)\}=0.\\ & \Rightarrow\; (2-\beta S)^*(2-\beta
S)=(2{\beta}^2P-\beta S)^*(2{\beta}^2P-\beta S).
\end{align*}
Since $r(S) < 2$, the operator $2-\beta S$ is invertible. Therefore,
we have $$ ((2-\beta S)^{-1})^*(2{\beta}^2P-\beta
S)^*(2{\beta}^2P-\beta S)(2-\beta S)^{-1}= I.$$ Therefore
$(2{\beta}^2P-\beta S)(2-\beta S)^{-1}$ and hence $(2{\beta}P-
S)(2-\beta S)^{-1}$ is an isometry for all $\beta\in\mathbb{T}$.

Conversely, let (5) hold. Then $(2{\beta}^2P-\beta S)(2-\beta
S)^{-1}$ is an isometry for every $\beta\in\mathbb{T}$.
Therefore,\begin{eqnarray*} & & ((2-\beta
S)^{-1})^*(2{\beta}^2P-\beta S)^*(2{\beta}^2P-\beta S)(2-\beta
S)^{-1}=I\\
& \mbox{ or } &  (2-\beta S)^*(2-\beta S)-(2{\beta}^2P-\beta
S)^*(2{\beta}^2P-\beta S)=0\\& \mbox{ or } & \rho(\beta
S,{\beta}^2P)=0\;,\quad \forall\;\beta\in\mathbb{T}.
\end{eqnarray*} Hence done.
\end{proof}

\begin{note} The $\Gamma$-isometry $(S_2, P_2)$ in the above proof
is a {\em pure} $\Gamma$-{\em isometry}.
    \end{note}
\begin{cor}
If $(S,P)$ is a $\Gamma$-isometry (respectively a $\Gamma$-unitary),
then $(rS,P)$ is also a $\Gamma$-isometry (respectively a
$\Gamma$-unitary) for $0\leq r \leq 1.$
\end{cor}

The following two results are remarkable in their simplicity to
characterize $\Gamma$-isometries.

\begin{lem} \label{goodlem}
A pair of bounded operators $(S,P)$ defined on $\mathcal{H}$ is a
pure $\Gamma$-isometry if and only if $S=C+C^*P$ for some pure
isometry $P$ and a bounded operator $C$ which commutes with $P$ and
$P^*$ and has numerical radius not greater than $1$.
\end{lem}
\begin{proof}
Let $(S,P)$ be a pure $\Gamma$-isometry. Then by Theorem 2.4 of
\cite{ay-jot}, $S$ and $P$ can be identified with $M_{\varphi}^E$
and $M_z^E$ respectively on $H^2(E)$ for some separable Hilbert
space $E$, where $\varphi(z)= G+G^*z$ for an operator $G$ defined on
$E$ such that $\omega(G)\leq1.$ Clearly $$
M_{\varphi}^E=M_{G+G^*z}^E=M_G^E+M_{G^*}^EM_z^E\equiv (I\otimes
G)+(I\otimes G^*)(M_z\otimes I) \;\mbox{on} \;
H^2(\mathbb{D})\otimes E\equiv H^2(E).$$ Therefore $S\equiv C+C^*P$
where $P=M_z\otimes I$ and $C=I\otimes G$. Obviously $P$ commutes
with $C, C^*$ and
$\omega(C)\leq1.$\\
Conversely, let $S=C+C^*P$ where $\omega(C)\leq1$ and $P$ is a pure
isometry which commutes with $C$ and $C^*$. Since $P$ is a pure
isometry, $P\equiv M_z^E$ on $H^2(E)$ and hence $C\equiv
M_{\varphi}^E$ on $H^2(E)$ for some $\varphi\in \mathcal{H}^{\infty}(\mathcal{B}(E))$, by the commutativity of $C$ and $P$.\\
Also since both of $M_{\varphi}^E$ and ${M_{\varphi}^E}^*$ commute
with $M_z^E$, the function $\varphi$ is a constant say equal to
$G_1$. Clearly $$M_{\varphi}^E\equiv (I\otimes G_1) \; \mbox{and}
\;M_z^E\equiv M_z\otimes I \;\mbox{on}\; H^2(\mathbb{D})\otimes E.$$
By the commutativity of $C$ and $P$ we have
$$S^*P=(C^*+P^*C)P=S.$$ Now
$$\omega(S)=\omega(C+C^*P)\leq
\omega(C)+\omega(C^*P)=\omega(I\otimes G_1)+\omega (M_z\otimes
G_1^*)\leq 1+ \omega(M_z\otimes G_1^*).$$ Since
$\omega(G_1^*)\leq1$, by Ando's result (Theorem \ref{ando}) there
exists a contraction $T$ such that
$$G_1^*=2(I-T^*T)^{1/2}T.$$ Considering the contraction $T_1=M_z\otimes
T$ we get \begin{align*} 2(I-T_1^*T_1)^{1/2}T_1 &=2(I\otimes
I-(M_z^*\otimes T^*)(M_z\otimes T))^{1/2}(M_z\otimes T)\\ &
=2(I\otimes I-I\otimes T^*T)^{1/2}(M_z\otimes T)\\& =2I\otimes
(I-T^*T)^{1/2}(M_z\otimes T)\\& = M_z\otimes
\{2(I-T^*T)^{1/2}T\}\\&=M_z\otimes G_1^*.
\end{align*} Therefore by Ando's result again, $\omega(M_z\otimes
G_1^*)\leq1.$ Thus we have $\omega(S)\leq2$. Therefore by Theorem
\ref{Gamma-isometry}-(3), $(S,P)$ is a $\Gamma$-isometry where $P$
is a pure isometry i.e, $(S,P)$ is a pure $\Gamma$-isometry.

\end{proof}
\begin{thm} \label{goodthm}
A pair of bounded operators $(S,P)$ defined on $\mathcal{H}$ is a
$\Gamma$-isometry if and only if $S=C+C^*P$ for some isometry $P$
and a bounded operator $C$ which commutes with $P$ and $P^*$ and has
numerical radius not greater than $1$.
\end{thm}
\begin{proof}
By Theorem \ref{Gamma-isometry}, $(S,P)$ is a $\Gamma$-isometry if
and only if $S=S_1\oplus S_2$ and $P=P_1 \oplus P_2$ where
$(S_1,P_1)$ and $(S_2,P_2)$ are $\Gamma$-unitary and pure
$\Gamma$-isometry
respectively.\\

Therefore $S_2=C+C^*P_2$ where $\omega(C)\leq1$ and $P_2$ is a pure
isometry which commutes with $C$ and $C^*$. Also by Theorem
\ref{G-unitary}, $S_1=U+U^*P_1$ where $U$ is a unitary
which commutes with $P_1$.\\
Choosing $C_1=U\oplus C$ we get $$S=S_1 \oplus S_2=
C_1+C_1^*(P_1\oplus P_2)=C_1+C_1^*P,$$ where $P$ commutes with
$C_1$, $C_1^*$ and obviously $\omega(C_1)\leq1$.
\end{proof}

\begin{obs}Let
$(S,P)$ be a $\Gamma$-contraction where $P$ is a projection. Then
$S$ and $P$ have the operator
matrices$$S\;=\;\begin{pmatrix}S_1&0\\0&S_2\end{pmatrix}\quad
P\;=\;\begin{pmatrix}I&0\\0&0\end{pmatrix}\;,$$ with respect to the
decomposition $\mathcal{H}=Ran(P)\oplus
Ker(P).$\end{obs}\vspace{0.12in}\begin{proof} Clearly $P$ has the
stated form as $P$ is a projection. Let $S=[S_{ij}]_{i,j=1}^2$ with
respect to the decomposition $\mathcal{H}=Ran(P)\oplus Ker(P).$ By
the commutativity of $S$ and $P$ it follows that $S_{12}=S_{21}=0.$
\end{proof}
\begin{obs}
If $(S,P)$ is a $\Gamma$-contraction where $P$ is a partial isometry
then $S-S^*P=\begin{pmatrix} 0&0\\{\ast}&S
\end{pmatrix}$ with respect to the decomposition
$\mathcal{H}=\overline{Ran}P^* \oplus Ker(P)$.
\end{obs}
\begin{proof}
Since $(S,P)$ is a $\Gamma$-contraction by Proposition \ref{prop1},
$\rho(\alpha S, {\alpha}^2P)\geq0$ for all $\alpha$ in $\mathbb{T}$
which implies that $$(I-P^*P)- \mbox{Re}\; \alpha(S-S^*P)\geq0.$$
Since $P$ is a partial isometry, $P^*P$ is a projection onto
$\overline{Ran}P^*=Ker(P)^{\bot}$. Therefore $I-P^*P$ is a
projection onto $Ker(P).$ So we have $P_{Ker(P)}-\mbox{Re}\; \alpha
(S-S^*P)\geq0$ for all $\alpha$ in $\mathbb{T}$. Therefore for $x\in
Ker(P)^{\bot}=\overline{Ran}P^*$ we have $P_{Ker(P)}(x)=0$ and hence
$$\mbox{Re}\;
\alpha(S-S^*P)|_{\overline{Ran}P^*}\leq0,\;  \mbox{for all}\; \alpha
\;\mbox{in}\; \mathbb{T}.$$ Therefore by Lemma \ref{basiclemma},
$(S-S^*P)|_{\overline{Ran}P^*}=0.$ Hence
$\overline{Ran}(S-S^*P)\subseteq Ker(P)$ and $S=\begin{pmatrix}
0&0\\ \ast &S
\end{pmatrix}$ with respect to the decomposition $\mathcal{H}= \overline{Ran}P^*\oplus
Ker(P)$.
\end{proof}
A canonical way of constructing a $\Gamma$-isometry is to consider
the Hardy space $H^2(\mathbb{D}^2)$ of the bidisc with the
reproducing kernel
$\frac{1}{(1-z_{1}\overline{w}_{1})(1-z_{2}\overline{w}_{2})}$. If
$M_{z_1}$ and $M_{z_2}$ are multiplications by the independent
variables $z_1$ and $z_2$ respectively, then $(M_{z_1} + M_{z_2} ,
M_{z_1} M_{z_2})$ is a $\Gamma$-isometry.

\section{$\Gamma$-contractions - examples}

Dilating a contraction operator to an isometry is well studied in
the history of dilation theory (see \cite{nazy}). For the class of
examples of $\Gamma$-contractions contained in this section, we
produce their $\Gamma$-isometric dilations.

\begin{defn} Let $(S,P)$ be a $\Gamma$-contraction on $\mathcal H$. A commuting pair of operators
$(T,V)$ acting on a Hilbert space $\mathcal{N}$ containing $
\mathcal{H}$ as a subspace is said to be a $\Gamma$-isometric
dilation of $(S,P)$ if $(T,V)$ is a $\Gamma$-isometry and
    $$T^*|_{\mathcal{H}}=S^*  \textup{ and } V^*|_{\mathcal{H}}=P^*.$$
   \end{defn} Thus $(T,V)$ is a $\Gamma$-isometric dilation of a $\Gamma$-contraction
     $(S,P)$ is same as saying that $(T^*,V^*)$ is a $\Gamma$-co-isometric extension of
     $(S^*,P^*)$. Moreover, the dilation will be called minimal if $$\mathcal N=\overline{\mbox{span}}\{V^nh: h\in \mathcal H\,\mbox{ and }
     n=0,1,2,\dots\}.$$ We shall see the existence and uniqueness
     of minimal $\Gamma$-isometric dilation in Theorem
     \ref{schaeffer}.

In this section we exhibit a new class of examples of
$\Gamma$-contractions and using a recent theorem of Douglas, Misra
and Sarkar (\cite{DMS}), find $\Gamma$-isometric dilations of some
of them. The main result of this section is Theorem \ref{example}.
\begin{lem} \label{fstlem}
Let $T_{1}$ and $T_{2}$ be to commuting contractions defined on
$\mathcal{H}$ and let $\mathcal{M}\subseteq\mathcal{H}$ be a
subspace invariant under $T_{1}+T_{2}$ and $T_{1}T_{2}$. Then
$((T_{1}+T_{2})|_{\mathcal{M}},T_{1}T_{2}|_{\mathcal{M}})$ is a
$\Gamma$-contraction.
\end{lem}

\begin{proof}
We have to show that $\Gamma$ is a spectral set for
$((T_{1}+T_{2})|_{\mathcal{M}},T_{1}T_{2}|_{\mathcal{M}})$, that is,
for any polynomial $p$ of two variables, \begin{align*}
\|p((T_{1}+T_{2})|_{\mathcal{M}},T_{1}T_{2}|_{\mathcal{M}})\|\leq\|p(z_{1},z_{2})\|_{{\infty},\Gamma}.
\end{align*}
Let $ \pi:\; \mathbb{C}^{2} \rightarrow \mathbb{C}^{2}$ be defined
as $$\pi(z_{1},z_{2})\; =\; (z_{1}+z_{2},z_{1}z_{2}).$$ Then by
\textit{von Neumann's} inequality in the bidisc $\mathbb{D}^{2}$, we
have,
\begin{align*}
\|p(\pi(T_{1},T_{2}))\|\leq \|p\circ\pi\|_{\infty,\mathbb{D}^2}\\
\mbox{or }\; \|p(T_{1}+T_{2},T_{1}T_{2})\|\leq \|p\|_{\infty,\Gamma}
\end{align*}
Certainly,
\begin{align*}
\|p((T_{1}+T_{2})|_{\mathcal{M}},T_{1}T_{2}|_{\mathcal{M}})\|\leq
\|p(T_{1}+T_{2},T_{1}T_{2})\|.
\end{align*}
Hence done.
\end{proof}\vspace{0.20in}
Let us see a particular example of this theorem. For $\lambda,\mu>1$
we define the weighted Bergman spaces
\begin{align} \label{wB}
\mathbb{A}^{(\lambda,\mu)}(\mathbb{D}^2)
&=\{f:\mathbb{D}^2\rightarrow\mathbb{C}\;:f \mbox{ is holomorphic
and }\\ \nonumber & \int_{\mathbb D^2
}|f(z_{1},z_{2})|^2(1-|z_{1}|^2)^{\lambda-2}(1-|z_{2}|^2)^{\mu-2}dm^{(\lambda,\mu)}(z_1,
z_2) < \infty\},
\end{align}
where  $m^{(\lambda,\mu)}$ is $\frac{(\lambda-1)(\mu-1)}{\pi^2}$
times the Lebesgue measure on $\mathbb D^2$. It is easy to verify
that $\mathbb{A}^{(\lambda,\mu)}(\mathbb D^2)$ is a Hilbert space.
For $f,g \in {\mathbb A}^{(\lambda,\mu)}(\mathbb D^2)$, define
 $$ \langle f,g \rangle_{\mathbb{A}^{(\lambda,\mu)}(\mathbb D^2)} =
  \int_{\mathbb D^2 }f(z_{1},z_{2})\overline{g(z_1,z_2)}(1-|z_{1}|^2)^{\lambda-2}(1-|z_{2}|^2)^
  {\mu-2}dm^{(\lambda, \mu)}(z_1 , z_2).$$

 Let $\Gamma^0$ denote the interior of $\Gamma$. Define the
 Hilbert space
\begin{align} \label{Gamma-spaces}
\mathbb{A}^{(\lambda,\mu)}({\Gamma}^0)
&=\{f:{\Gamma}^0\rightarrow\mathbb{C}: f \mbox{ is holomorphic and }
(f\circ\pi) \det
J_{\pi}\in\mathbb{A}^{(\lambda,\mu)}(\mathbb{D}^2)\},
\end{align}
with $$\langle f,g \rangle_{\mathbb{A}^{(\lambda,\mu)}(\Gamma^0)}=
\langle (f\circ\pi)\det J_{\pi},(g\circ \pi)\det J_{\pi}
\rangle_{\mathbb {A}^{(\lambda,\mu)}(\mathbb D^2)}\,,$$

where $J_{\pi}=\begin{pmatrix} 1&1\\
z_{2}&z_{1}\\
\end{pmatrix}$ is the Jacobian of the map $\pi(z_{1},z_{2})= (z_{1}+z_{2},z_{1}z_{2})$
so that $\det J_{\pi}=(z_{1}-z_{2})$. Let $M_s^{(\lambda, \mu)}$ and
$M_p^{(\lambda, \mu)}$ be the multiplication operators on
$\mathbb{A}^{(\lambda,\mu)}({\Gamma}^0)$ by the co-ordinate
functions $s$ and $p$, respectively, where $(s,p) \in \Gamma^0$.
 For $\lambda = \mu,$ we denote
$\mathbb{A}^{(\lambda,\lambda)}(\mathbb{D}^2)$ by
$\mathbb{A}^{(\lambda)}(\mathbb{D}^2)$,
$\mathbb{A}^{(\lambda,\lambda)}({\Gamma}^0)$ by
$\mathbb{A}^{(\lambda)}({\Gamma}^0)$, $M_s^{(\lambda, \lambda)}$ by
$M_s^{(\lambda)}$ and $M_p^{(\lambda, \lambda)}$ by
$M_p^{(\lambda)}$. The following lemma serves the purpose of showing
that the operator pair $(M_s^{(\lambda,\mu)},M_p^{(\lambda,\mu)})$,
which is obviously a commuting pair, is a $\Gamma$-contraction.
\begin{lem}\label{key}
For integers $m,n \geq 0$, let $\widetilde{e_{mn}}$ and
$\widetilde{f_{mn}}$ be the functions defined on $\mathbb{D}^2$ by
$$
\widetilde{e_{mn}}(z_{1},z_{2})= z_{1}^mz_{2}^n-z_{1}^nz_{2}^m
\mbox{~~and~~} \widetilde{f_{mn}}(z_{1},z_{2})=
z_{1}^mz_{2}^n+z_{1}^nz_{2}^m.
$$
Let  $\mathbb{A}_\mathrm{a}^{(\lambda,\mu)} (
\mathbb{D}^2):=\overline{\mathrm{span}}\{\widetilde{e_{mn}}: m>n\geq
0\} $ and $\mathbb{A}_\mathrm{s}^{(\lambda,\mu)} (
\mathbb{D}^2):=\overline{\mathrm{span}}\{\widetilde{f_{mn}}: m\geq
n\geq 0\} $ be  subspaces of $\mathbb{A}^{(\lambda,\mu)} (\mathbb
D^2)$. Then
\begin{itemize}
\item[(1)] both $\mathbb{A}_\mathrm{a}^{(\lambda,\mu)} (
    \mathbb{D}^2) $ and $\mathbb{A}_\mathrm{s}^{(\lambda,\mu)}
    ( \mathbb{D}^2) $ are
     invariant subspaces of $\mathbb{A}^{(\lambda,\mu)} (
    \mathbb{D}^2) $ under $M_{z_1+z_2}$ and $M_{z_1z_2}$\;;
\item[(2)] the restrictions of the pair
    $(M_{z_1+z_2},M_{z_1z_2})$ to the invariant subspaces
    $\mathbb{A}_\mathrm{a}^{(\lambda,\mu)}(\mathbb D^2)$ and
    $\mathbb{A}_\mathrm{s}^{(\lambda,\mu)}(\mathbb D^2)$  are
    $\Gamma$-contractions, call them
    $(S_a^{(\lambda,\mu)},P_a^{(\lambda,\mu)})$ and
    $(S_s^{(\lambda,\mu)} P_s^{(\lambda,\mu)})$ respectively.
    As usual, for $\lambda=\mu$ we use just one index.
\item[(3)] there is an isometry $U$ from
    $\mathbb{A}^{(\lambda,\mu)}({\Gamma}^0)$ onto
    $\mathbb{A}_\mathrm{a}^{(\lambda,\mu)}(\mathbb{D}^2)$ such
    that $UM_{s}^{(\lambda,\mu)}U^*= S_a^{(\lambda,\mu)}$ and
    $UM_{p}^{(\lambda,\mu)}U^*= P_a^{(\lambda,\mu)}$.
\end{itemize}
\end{lem}\vspace{0.07in}
\begin{proof}
 Note that
\begin{align*}
(z_{1}+z_{2})(z_{1}^mz_{2}^n-z_{1}^nz_{2}^m) &=
z_{1}^{m+1}z_{2}^n-z_{1}^{n+1}z_{2}^m+z_{1}^mz_{2}^{n+1}-z_{1}^nz_{2}^{m+1}\\
&=
(z_{1}^{m+1}z_{2}^n-z_{1}^nz_{2}^{m+1})+(z_{1}^mz_{2}^{n+1}-z_{1}^{n+1}z_{2}^m)
\end{align*}
Again
$$z_{1}z_{2}(z_{1}^mz_{2}^n-z_{1}^nz_{2}^m)=(z_{1}^{m+1}z_{2}^{n+1}-z_{1}^{n+1}z_{2}^{m+1}).$$\\

So $\mathbb{A}_\mathrm{a}^{(\lambda,\mu)}(\mathbb D^2)$ is invariant
under both of the multiplication operators $ M_{z_{1}+z_{2}}$ and
$M_{z_{1}z_{2}}$. We can show similarly that
$\mathbb{A}_\mathrm{s}^{(\lambda,\mu)}(\mathbb D^2)$ is invariant
under both the operators $ M_{z_{1}+z_{2}}$ and $M_{z_{1}z_{2}}$.
Hence by Lemma \ref{fstlem}, (1) and (2) above are proved. To prove
(3), define $$U:
\mathbb{A}^{(\lambda,\mu)}({\Gamma}^0)\longrightarrow
\mathbb{A}^{(\lambda,\mu)}(\mathbb{D}^2)$$ by
\begin{align*}
Uf=(f\circ\pi)\det J_{\pi}.
\end{align*}
That $U$ is an isometry follows from the definitions of norms on the
corresponding spaces. It is easy to check by direct computation that
$U$ intertwines $M_{s}^{(\lambda, \mu)}$ with $S_a^{(\lambda,\mu)}$
and $M_{p}^{(\lambda, \mu)}$ with $P_a^{(\lambda,\mu)}$.
\end{proof}
\begin{rem}\rm
We observe that ${\langle{\widetilde {e_{mn}}},{\widetilde
{f_{mn}}}\rangle}_{\mathbb A^{(\lambda,\mu)}(\mathbb
D^2)}=\frac{m!n!}{(\lambda)_m(\mu)_n}-\frac{m!n!}{(\mu)_m(\lambda)_n},$
for $m>n\geq 0,$ where
$(\lambda)_m=\frac{\lambda(\lambda-1)(\lambda-2)\dots(\lambda-m+1)}{m!}.$
 Therefore the subspaces $\mathbb
A^{(\lambda,\mu)}_\mathrm{a}(\mathbb D^2)$ and $\mathbb
A^{(\lambda,\mu)}_\mathrm{s}(\mathbb D^2)$ of $\mathbb
A^{(\lambda,\mu)}(\mathbb D^2)$ are mutually orthogonal if and only
if $\lambda=\mu.$
\end{rem}

Consider the weighted Bergman space $\mathbb A^{(\lambda)}(\mathbb
D^2)$, as defined in (\ref{wB}), on the bidisc for $\lambda> 1$ and
its subspaces
$$\mathbb A_{\mathrm{a}}^{(\lambda)}(\mathbb D^2):=\overline{\mathrm{span}}\{z_1^mz_2^n-z_1^nz_2^m:
m>n\geq 0, (z_1,z_2)\in \mathbb D^2\}$$ and
$$\mathbb A_{\mathrm{s}}^{(\lambda)}(\mathbb D^2):=\overline{\mathrm{span}}\{z_1^mz_2^n+z_1^nz_2^m: m\geq n\geq 0, (z_1,z_2)\in \mathbb D^2\}.$$
They are mutually orthogonal and $\mathbb A^{(\lambda)}(\mathbb
D^2)=\mathbb A_{\mathrm{a}}^{(\lambda)}(\mathbb D^2)\oplus \mathbb
A_{\mathrm{s}}^{(\lambda)}(\mathbb D^2).$ Let

\begin{equation} \label{defnL}
f_{mn}(z_1,z_2)=\left\{%
\begin{array}{ll} \sqrt{\frac{(\lambda)_m(\lambda)_n}{2 (m!n!)}}(z_1^mz_2^n+z_1^nz_2^m)
     & {\mathrm{for}~ m>n\geq 0;} \\
    \sqrt{\frac{(\lambda)_n}{n!}}(z_1z_2)^n & {\mathrm{for}~ m=n\geq 0.} \\
\end{array} \right.
\end{equation}
Clearly, $\{f_{mn}\}_{m\geq n\geq 0}$ is an orthonormal basis for
the Hilbert space $\mathbb A_{\mathrm{s}}^{(\lambda)}(\mathbb D^2)$.

\begin{prop}
 The Hilbert space ${\mathbb A}_{\mathrm{s}}^{(\lambda)}(\mathbb D^2)$ is a
 reproducing kernel Hilbert space with its reproducing kernel $K^{(\lambda)}_{\mathrm{s}}$  given by the formula:
\begin{eqnarray} \label{formulaforklambdasym}
 K^{(\lambda)}_{\mathrm{s}} (\boldsymbol z, \boldsymbol w) =
 \frac{1}{2} (1-z_1\bar w_1)^{-\lambda} (1-z_1\bar w_2)^{-\lambda} + \frac{1}{2} (1-z_1\bar w_2)^{-\lambda} (1-z_2 \bar w_1)^{-\lambda},
 \end{eqnarray}
 where $ \boldsymbol z=(z_1,z_2)$ and $\boldsymbol w=(w_1,w_2)$
 are in $\mathbb {D}^2$
\end{prop}
\begin{proof}
We shall prove this by expanding the right hand side of the formula
(\ref{formulaforklambdasym}) in terms of the basis elements
$f_{mn}$.
 For $\boldsymbol z, \boldsymbol w\in \mathbb D^2,$ we have
\begin{eqnarray*} \lefteqn{K^{(\lambda)}_{\mathrm{s}}(\boldsymbol z,\boldsymbol w)
= \displaystyle\sum_{m\geq n\geq 0} f_{mn}(z_1,z_2)\overline{f_{mn}(w_1,w_2)}} \\
&=&  \displaystyle\sum_{m> n\geq 0}
f_{mn}(z_1,z_2)\overline{f_{mn}(w_1,w_2)}
+  \displaystyle\sum_{ n\geq 0} f_{nn}(z_1,z_2)\overline{f_{nn}(w_1,w_2)} \\
&=& \frac{1}{2} \displaystyle\sum_{\stackrel{m,n\geq 0}{m\neq n}}
f_{mn}(z_1,z_2)\overline{f_{mn}(w_1,w_2)}
+ \displaystyle\sum_{ n\geq 0} f_{nn}(z_1,z_2)\overline{f_{nn}(w_1,w_2)} \\
&=& \frac{1}{4}\displaystyle\sum_{\stackrel{m,n\geq 0}{m\neq
n}}{\frac{(\lambda)_m(\lambda)_n}{m!n!}}
\Big((z_1\bar w_1)^m(z_2\bar w_2)^n+ (z_1\bar w_2)^m(z_2\bar w_1)^n+(z_2\bar w_1)^m(z_1\bar w_2)^n \\
& &\hspace{2 cm}+(z_2\bar w_2)^m(z_1\bar w_1)^n\Big)+
\displaystyle\sum_{n\geq 0}
\frac{(\lambda)^2_n}{(n!)^2}(z_1\bar w_1)^n(z_2\bar w_2)^n \\
&=& \frac{1}{2}(1-z_1\bar w_1)^{-\lambda}(1-z_1\bar
w_2)^{-\lambda}+\frac{1}{2}(1-z_1\bar w_2)^{-\lambda}(1-z_2\bar
w_1)^{-\lambda}
\end{eqnarray*}
\end{proof}

Edigarian and Zwonek found the Bergman kernel for symmetrized
polydisc, see \cite{edi-zwo}. We shall need explicit formulae for
the reproducing kernels of the weighted Bergman spaces $\mathbb
A^{(\lambda)}(\Gamma^0)$, as defined in (\ref{Gamma-spaces}). These
have been extensively studied in \cite{MSZ}. We recall only some
relevant facts here. For $\lambda>1$, the reproducing kernel for the
weighted Bergman space $\mathbb A^{(\lambda)}(\Gamma^0)$ on the
interior of the symmetrized bidisc $\Gamma^0$ is given by
\begin{equation}\label{Berg}
 \mathbf B_{\Gamma^0}^{(\lambda)}(\pi(\boldsymbol z),\pi(\boldsymbol w))
 ={\frac{1}{\lambda}}\frac{\{(1-z_1\bar w_1)^{-\lambda}(1-z_1\bar w_2)^{-\lambda}-(1-z_1\bar w_2)^{-\lambda}(1-z_2\bar w_1)^{-\lambda}\}}
{(z_1-z_2)(\bar w_1-\bar w_2)},\quad \boldsymbol z, \boldsymbol w\in
\mathbb D^2.
\end{equation}

The kernel above remains a positive definite kernel for $\lambda=1$.
This prompted the authors of \cite{MSZ} to define the Hardy space
$H^2(\Gamma^0)$ of the symmetrized bidisc to be the reproducing
kernel Hilbert space whose kernel is
\begin{equation}\label{Sze}
 \mathbb S_{\Gamma^0}(\pi(\boldsymbol z),\pi(\boldsymbol w))
 =\frac{(1-z_1\bar w_1)^{-1}(1-z_1\bar w_2)^{-1}-(1-z_1\bar w_2)^{-1}(1-z_2\bar w_1)^{-1}}
 {(z_1-z_2)(\bar w_1-\bar w_2)}, \quad \boldsymbol z, \boldsymbol w\in \mathbb D^2.
\end{equation}

\begin{lem}\label{quotient}
 The ratio $\mathbb S_{\Gamma^0}^{-1}\mathbf B_{\Gamma^0}^{(n)}$ of the reproducing kernel of the weighted Bergman space
 $\mathbb A^{(n)}(\Gamma^0)$ with the reproducing kernel of the Hardy space $H^2(\Gamma^0)$ is a positive definite kernel
 for all positive integers $n$.

\end{lem}
\begin{proof}
For $\boldsymbol z,\boldsymbol w\in \mathbb D^2$, from \eqref{Berg}
and \eqref{Sze}, we have

\begin{align*}
\mathbb S_{\Gamma^0}^{-1}\mathbf B_{\Gamma^0}^{(n)}(\pi(\boldsymbol
z),\pi(\boldsymbol w)) &=\frac{1}{\lambda}\frac{(1-z_1\bar
w_1)^{-n}(1-z_2\bar w_2)^{-n}
-(1-z_1\bar w_2)^{-n}(1-z_2\bar w_1)^{-n}}{(1-z_1\bar w_1)^{-1}(1-z_2\bar w_2)^{-1}-(1-z_1\bar w_2)^{-1}(1-z_2\bar w_1)^{-1}}\\
&=\frac{1}{\lambda}\displaystyle\sum_{k=0}^{n-1} a^{n-1-k}b^k,
\end{align*} where $a=(1-z_1\bar w_1)^{-1}(1-z_2\bar w_2)^{-1}$ and
$b=(1-z_1\bar w_2)^{-1}(1-z_2\bar w_1)^{-1}.$  Clearly, the last
expression can be expressed as a polynomial in $ab$ and $a^k+b^k$
for $k=1,\ldots,n-1.$ Since $ab=\mathbb S_{\Gamma^0}(\pi(\boldsymbol
z),\pi(\boldsymbol w))$ is the reproducing kernel for the Hilbert
spaces $H^2(\Gamma^0)$ and $a^k+b^k=K^{(k)}_{\mathrm{s}}(\boldsymbol
z,\boldsymbol w)$ is the reproducing kernel for the Hilbert space
$\mathbb A^{(k)}_{\mathrm{s}}(\mathbb D^2),$ they both are positive
definite. Recalling that pointwise product and sum of two positive
definite kernels are again positive definite kernels we conclude
that $\mathbb S_{\Gamma^0}^{-1}\mathbf B_{\Gamma^0}^{(n)}$ is a
positive definite kernel.
\end{proof}

By $H^2(\mathbb{D}^2)$, we shall denote the Hardy space of the
bidisc. For convenience of notation, we shall also call it
$\mathbb{A}^{(1)}(\mathbb{D}^2)$. This will enable us to talk about
the operator pairs $(S_a^{(1)}, P_a^{(1)})$ and $(S_s^{(1)},
P_s^{(1)})$.
\begin{lem}
 The pair $(M_s^{H}, M_p^{H})$ of multiplication operators  on $H^2(\Gamma^0)$ by the co-ordinate functions is a $\Gamma$-isometry.
\end{lem}
\begin{proof}
Let  $H^2_{\mathrm{a}}(\mathbb
D^2):=\overline{\mathrm{span}}\{z_1^mz_2^n-z_1^nz_2^m: m\geq n\geq
0, (z_1,z_2)\in \mathbb D^2\}$
 and  $H^2_{\mathrm{s}}(\mathbb D^2):=\overline{\mathrm{span}}\{z_1^mz_2^n+z_1^nz_2^m: m\geq n\geq 0, (z_1,z_2)\in \mathbb D^2\}.$ Clearly,
 $H^2(\mathbb D^2)=H^2_{\mathrm{a}}(\mathbb D^2)\oplus H^2_{\mathrm{s}}(\mathbb D^2).$ For $\lambda=\mu=1,$ analogous arguments as in
Lemma \ref{key}, shows that
\begin{enumerate}
 \item[(i)] the subspaces $ H^2_\mathrm{a}(\mathbb{D}^2)$ and
     $ H^2_\mathrm{s}(\mathbb{D}^2)$  are invariant subspaces
     of $ H^2(\mathbb{D}^2)$ under $M_{z_1+z_2}$ and
     $M_{z_1z_2};$ \item[(ii)] there is an isometry $U$ from
    $H^2({\Gamma}^0)$ onto $H^2_\mathrm{a}(\mathbb{D}^2)$ such
    that $UM_{s}^HU^*= S_a^{(1)}$ and $UM_{p}^HU^*=
    P_a^{(1)}$.
\end{enumerate}
By Theorem \ref{G-unitary}$-(2)$, the pair
$(M_{z_1+z_2},M_{z_1z_2})$ is a $\Gamma$-unitary on $L^2( \mathbb
T^2).$ Moreover,
$S_a^{(1)}=M_{z_1+z_2}|_{H^2_\mathrm{a}(\mathbb{D}^2)}$ and
$P_a^{(1)}=M_{z_1z_2}|_{H^2_\mathrm{a}(\mathbb{D}^2)}.$ So
$(S_a^{(1)}, P_a^{(1)})$ is a $\Gamma$-isometry. Noting that $U$ is
a unitary, it follows from (ii) that $(M_s^H, M_p^H)$ is a
$\Gamma$-isometry.
\end{proof}

Recall that $(M^{(\lambda)}_s, M^{(\lambda)}_p)$ denotes the
commuting pair of multiplication operators by the coordinate
functions $s$ and $ p$, respectively, on the Hilbert space $\mathbb
A^{(\lambda)}(\Gamma^0)$ for $\lambda>1$. For $\lambda=1$, this
space is the Hardy space $H^2(\Gamma^0)$ and the operator pair is
$(M_s^H,M_p^H)$. Thus, by Lemma \ref{quotient} and Theorem 6 of
\cite{DMS}, we have proved the following theorem which is the main
result of this section.
\begin{thm}\label{example}
 For every positive integer $n,$ the $\Gamma$-contraction
 $(M_s^{(n)},M_p^{(n)})$acting on $\mathbb A^{(n)}(\Gamma^0)$ can
 be dilated to the $\Gamma$-isometry $(M_s^H\otimes I_{\mathcal L},M_p^H\otimes I_{\mathcal
 L})$ on $H^2(\Gamma^0)\otimes \mathcal L$ for some Hilbert space
 $\mathcal L$.
\end{thm}
\noindent Recalling the notations from Lemma \ref{key}, we have the
following corollary.
\begin{cor}
For every positive integer $n,$ the commuting pair of operators
$(S_a^{(n)}, P_a^{(n)})$ acting on the Hilbert space $\mathbb
A^{(n)}_\mathrm{a}(\mathbb D^2)$ has a $\Gamma$-isometric dilation
to the commuting pair of operators $(S_a^{(1)}, P_a^{(1)})$ on the
Hilbert space $H^2_\mathrm{a}(\mathbb D^2)\otimes \mathcal{L}$ for
some Hilbert space $\mathcal L.$
\end{cor}
\begin{proof}
Observing that the isometry $U$ in part 3 of Lemma \ref{key} is
actually a unitary the proof follows from Theorem \ref{example}.
\end{proof}

\begin{lem}
The $\Gamma$-isometric dilation $(S_a^{(1)},P_a^{(1)})$ on the
Hilbert space $H^2_\mathrm{a}(\mathbb D^2)\otimes \mathcal{L}$ of
the commuting pair of operators $(S_a^{(n)},P_a^{(n)})$ on the
Hilbert space $\mathbb A^{(n)}_\mathrm{a}(\mathbb D^2)$ is minimal.
\end{lem}
\begin{proof}
 To prove minimality, we need to show that $\overline{\mathrm{span}}
 \{P_1^kh: h\in H^2_\mathrm{a}(\mathbb D^2), k\geq 0\}=H^2_\mathrm{a}(\mathbb D^2).$
Recalling that $\widetilde {e_{mn}}(z_1,z_2)=z_1^mz_2^n-z_1^nz_2^m$
and $ H^2_\mathrm{a}(\mathbb
D^2)=\overline{\mathrm{span}}\{\widetilde {e_{mn}}: m>n\geq 0\}$, it
suffices to show that $\overline{\mathrm{span}}\{P_1^k(\widetilde
{e_{mn}}): m>n\geq 0,k\geq 0\} =\overline{\mathrm{span}}\{\widetilde
{e_{mn}}: m>n\geq 0\}.$ Since $P_1^k(\widetilde {e_{mn}})=\widetilde
{e_{m+k,n+k}}$, we have $\overline{\mathrm{span}}\{P_1^k(\widetilde
{e_{mn}}): m>n\geq 0,k\geq 0\} =\overline{\mathrm{span}}\{\widetilde
{e_{m+k,n+k}}: m>n\geq 0,k\geq
0\}=\overline{\mathrm{span}}\{\widetilde {e_{mn}}: m>n\geq 0\}.$
Hence the proof is complete.
\end{proof}

We have a corollary of the above lemma.

\begin{cor}
The dilation $(M_s^H\otimes I_{\mathcal L},M_p^H\otimes I_{\mathcal
L})$ on the Hilbert space $H^2(\Gamma^0)\otimes \mathcal{L}$ of the
commuting pair of operators $(M_s^{(n)},M_p^{(n)})$ on the Hilbert
space $\mathbb A^{(n)}(\Gamma^0)$ is minimal.
\end{cor}
\begin{proof}
 Set $\epsilon_{mn}(s,p)=\epsilon\circ\pi(z_1,z_2)=\frac{z_1^mz_2^n-z_1^nz_2^m}{z_1-z_2}$
 for $(s,p)\in \Gamma^0, (z_1,z_2)\in \mathbb D^2.$ So
$H^2(\Gamma^0)=\overline{\mathrm{span}}\{\epsilon_{mn}:m>n\geq 0\}$
and $M_p\epsilon_{mn}=\epsilon_{m+1,n+1}.$ Now analogous arguments
as in the previous corollary shows that
$$\overline{\mathrm{span}}\{M_p^k\epsilon_{mn}:m>n\geq 0, k\geq 0\}
=\overline{\mathrm{span}}\{\epsilon_{mn}:m>n\geq
0,\}=H^2(\Gamma^0).$$ This proves minimality of $(M_s^H\otimes
I_{\mathcal L}, M_p^H \otimes I_{\mathcal L})$ on the Hilbert space
$H^2(\Gamma^0)\otimes\mathcal L.$
\end{proof}

We move on to general discussion of dilation in the next section.

\section{Dilation}
As in many occasions in operator theory, in our case too, finding a
solution to an operator equation turns out to be of utmost
importance. As is clear by now, a crucial role in deciphering the
structure of a $\Gamma$-contraction $(S,P)$ is played by the
operator $S-S^*P$. For a pair $(S,P)$ of commuting bounded operators
with $\|P\|\leq1$, we shall denote from now on by $\Sigma$ and
$\Sigma_*$, the operators $S-S^*P$ and $S^*-SP^*$ respectively. We
denote by $\mathbf{D}_P$ and $\mathcal{D}_P$ the operator
$(I-P^*P)^{\frac{1}{2}}$ and its range closure respectively. For a
pair of commuting bounded operators $(S,P)$ with $\|P\|\leq1$, the
{\em fundamental equation} is defined to be
\begin{eqnarray}\label{e0} \Sigma=\mathbf{D}_PX\mathbf{D}_P,\;\mbox{where }X\in\mathcal B(\mathcal{D}_P)\end{eqnarray}
and the same for the pair $(S^*,P^*)$ is \begin{eqnarray} \label{e1}
{\Sigma}_*=\mathbf{D}_{P^*}Y\mathbf{D}_{P^*},\; \mbox{where
}Y\in\mathcal B(\mathcal{D}_{P^*}).\end{eqnarray} We start with a
pivotal theorem which guarantees the existence and uniqueness of
solutions of such equations for $\Gamma$-contractions. The proof of
the theorem needs the following lemma and the proof of the lemma
given here is from a private communication with Michael A.
Dritschel.

\begin{lem} \label{imp-lemma}
Let $\Sigma$ and $D$ be two bounded operators on $\mathcal H$. Then
$$ DD^*\geq \textup{ Re }(e^{i\theta} \Sigma) \quad \textit{ for all
} \theta \in [0,2\pi)
$$
if and only if there is $F \in \mathcal B(\mathcal D_*)$ with
numerical radius of $F$ not greater than one such that $\Sigma
=DFD^*$, where $\mathcal D_*=\overline{\textup{Ran }}D^* $.
\end{lem}

The proof of this result needs the operator Fejer-Riesz
factorization theorem (Theorem 2.1 of \cite{DW}) along with
Douglas's lemma (Lemma 2.1 of \cite{DMP}) and the familiar result
that an operator $X$ has numerical radius not greater than one if
and only if $\textup{Re }\beta X \leq I$ for all complex numbers
$\beta$ of modulus 1 (Lemma \ref{basicnrlemma}).

\textit{Proof of Lemma} \ref{imp-lemma}. Let there be an operator
$F\in\mathcal B(\mathcal H)$ with numerical radius not bigger than
one such that $\Sigma =DFD^*$. Since $I-\textup{Re
}(e^{i\theta}F)\geq 0$, for all $\theta\in [0,2\pi),$ we have
$$ D(I-\textup{Re }e^{i\theta}F)D^*\geq 0,\quad
\textup{for all } \theta.$$ So we have
$$ DD^*\geq D\textup{Re }(e^{i\theta}F)D^* =\textup{Re }(e^{i\theta} \; DFD^*)=\textup{Re }(e^{i\theta}\Sigma)$$
for all $\theta \in [0,2\pi)$.

The nontrivial part of this lemma, however, is the converse of the
above. Suppose that $DD^*\geq \textup{Re }(e^{i\theta} \Sigma )$ for
all $\theta\in [0,2\pi)$. This means that the Laurent polynomial
$$  DD^* - \frac{1}{2} ( z \Sigma + \bar{z} \Sigma^* )$$
is non-negative for $z$ on the unit circle. By the operator
Fej\'{e}r-Riesz theorem (Theorem 2.1 of \cite{DW}) we thus have a
factorization
$$ DD^*- \frac{1}{2} (z \Sigma + \bar z \Sigma^*) = (X-zY)(X^*-\bar zY^*),\quad
|z|=1,$$ with $X,Y\in \mathcal B(\mathcal H).$ Thus $DD^*=XX^*+YY^*$
and $\Sigma = 2YX^*$. Since $DD^* \ge XX^*$ and $DD^* \ge YY^*$,
Douglas's lemma tells us that there exist contractions $Q$ and $R$
such that $X=DQ$ and $ Y=DR$. Thus $\Sigma = DFD^*$ for
$$ F= P_{\mathcal D_*} 2RQ^*|_{\mathcal D_*}.$$

To show that the numerical radius of $F$ is not greater than one,
note that $ DD^*\geq \textup{Re }(e^{i\theta} \Sigma )=\textup{Re
}(e^{i\theta}DFD^*) $ for all $\theta\in [0,2\pi)$ which implies
that
$$ D(I_{\mathcal D_*} - \textup{Re }(e^{i\theta} F))D^*\geq 0, \quad \textup{for all }\theta\in
[0,2\pi).$$ Hence
$$ \langle ( I_{\mathcal D_*} - \textup{Re }( e^{i\theta} F )) D^* h , D^*h \rangle =
\langle D (I_{\mathcal D_*} - \textup{Re }( e^{i\theta} F )) D^* h
, h \rangle \geq 0$$ for all $\theta\in [0,2\pi)$ and as a
consequence, the numerical radius of $F$ is no bigger than one.
\qed

Now here is the theorem which guarantees the existence and
uniqueness of solution to the fundamental equation of a
$\Gamma$-contraction.

\begin{thm}\label{greatlem}
Let $(S,P)$ be a $\Gamma$-contraction. Then there is a unique
solution $A$ to its fundamental equation
$$S-S^*P=\mathbf{D}_PX\mathbf{D}_P.$$ Moreover, $A$ has
numerical radius less than or equal to one.
\end{thm}
\begin{proof}
Since $(S,P)$ is a $\Gamma$-contraction, by Proposition \ref{prop1},
we have
$$\rho(\alpha S,{\alpha}^2P)\geq 0 \quad \textup{for all } \alpha\in \overline{\mathbb D}. $$
So in particular for all $\beta$ with modulus 1, we have $
\mathbf{D}_P^2-\textup{Re }\beta(S-S^*P)\geq 0$. Therefore by Lemma
\ref{imp-lemma}, there exists an operator $A\in \mathcal B(\mathcal
D_P)$ with numerical radius not greater than one such that
$S-S^*P=\mathbf{D}_PA\mathbf{D}_P$.

For uniqueness let there be two such solutions $A_1$ and $A_2$. Then
$$\mathbf{D}_P\tilde{A}\mathbf{D}_P=0, \quad \textup{where }
\tilde{A}=A_1-A_2 \in \mathcal B(\mathcal D_P).$$ Then $$\langle
\tilde{A}\mathbf{D}_Ph,\mathbf{D}_Ph^{\prime} \rangle=\langle
\mathbf{D}_P\tilde{A}\mathbf{D}_Ph,h^{\prime} \rangle =0 $$ which
shows that $\tilde{A}=0$ and hence $A_1=A_2$.
\end{proof}

This theorem allows us to construct an explicit $\Gamma$-isometric
dilation of a $\Gamma$-contraction, which is one of our main results
and is shown in the following theorem.

\begin{thm} \label {schaeffer}
Let $(S,P)$ be a $\Gamma$-contraction on a Hilbert space
$\mathcal{H}$. Let $A$ be the unique solution of the fundamental
equation (\mbox{\ref{e0}}) and let
$\mathcal{K}_0=\mathcal{H}\oplus\mathcal{D}_{p}\oplus\mathcal{D}_{p}\oplus\mathcal{D}_{p}\oplus\dots=\mathcal{H}\oplus
l^2(\mathcal{D}_{p})$. Consider the operators $T_A,V_0$ defined on
$\mathcal{K}_0$ by \begin{align*} &
T_A(h_0,h_1,h_2,\dots)=(Sh_0,A^*\mathbf{D}_Ph_0+Ah_1,A^*h_1+Ah_2,A^*h_2+Ah_3,\dots)\\
& V_0(h_0,h_1,h_2,\dots)=(Ph_0,\mathbf{D}_Ph_0,h_1,h_2,\dots).
\end{align*} Then \begin{enumerate}
\item $(T_A,V_0)$ is a $\Gamma$-isometric dilation of $(S,P)$.
\item If $(\widehat{T},V_0)$ on $\mathcal{K}_0$ is a
    $\Gamma$-isometric dilation of $(S,P)$, then
    $\widehat{T}=T_A$.
\item If $(T,V)$ is a $\Gamma$-isometric dilation of $(S,P)$ where
    $V$ is a minimal isometric dilation of $P$, then $(T,V)$ is
    unitarily equivalent to $(T_A,V_0)$.
\end{enumerate} Thus $(2)$ and $(3)$ guarantee the uniqueness of
$\Gamma$-isometric dilation $(T,V)$ of $(S,P)$ where $V$ is minimal
isometric dilation of $P$.
\end{thm}

\begin{proof}

\textbf{(1)} It is evident from the definition that $V_0$ on
$\mathcal{K}_0$ is the minimal isometric dilation of $P$. Obviously
$T_A^*$ and $V_0^*$ are defined on $\mathcal{K}_0$ as
\begin{align*}
&
T_A^*(h_0,h_1,h_2,\dots)=(S^*h_0+\mathbf{D}_PAh_1,A^*h_1+Ah_2,A^*h_2+Ah_3,\dots)\\
& V_0^*(h_0,h_1,h_2,\dots)=(P^*h_0+\mathbf{D}_Ph_1,h_2,h_3,\dots).
\end{align*}
The space $\mathcal{H}$ can be embedded inside $\mathcal{K}_0$ by
the map $h\mapsto (h,0,0,\dots)$. It is clear that $\mathcal{H}$,
considered as a subspace of $\mathcal{K}_0$ is co-invariant under
$T_A$ and $V_0$ and
$T_A^*|_{\mathcal{H}}=S^*, V_0^*|_{\mathcal{H}}=P^*$.\\
Since $V_0$ is an isometry, in order to show that $(T_A,V_0)$ is a
$\Gamma$-{\em isometric dilation} of $(S,P)$ one has to justify (by
virtue of Theorem \ref{Gamma-isometry}-(3)) the following:
\begin{itemize}
\item[(a)]$T_AV_0=V_0T_A$\item[(b)]
    $T_A=T_A^*V_0$\item[(c)]$r(T_A) \leq2$.
\end{itemize}
\begin{align*}
T_AV_0(h_0,h_1,h_2,\dots)&=T_A(Ph_0,\mathbf{D}_Ph_0,h_1,h_2,\dots)\\
&=(SPh_0,A^*\mathbf{D}_Ph_0+A\mathbf{D}_Ph_0,A^*\mathbf{D}_Ph_0+Ah_1,A^*h_1+Ah_2,A^*h_2+Ah_3,\dots).
\end{align*}
\begin{align*}
V_0T_A(h_0,h_1,h_2,\dots)&=V_0(Sh_0,A^*\mathbf{D}_Ph_0+Ah_1,A^*h_1+Ah_2,A^*h_2+Ah_3,\dots)\\
&=(PSh_0,\mathbf{D}_PSh_0,A^*\mathbf{D}_Ph_0+Ah_1,A^*h_1+Ah_2,A^*h_2+Ah_3,\dots).
\end{align*}
Let $G=A^*\mathbf{D}_PP+A\mathbf{D}_p-\mathbf{D}_pS $. Then $G$ is
defined from $\mathcal{H}\rightarrow\mathcal{D}_P$. Since $A$ is a
solution of the equation (\ref{e0}), we have
\begin{align*}\mathbf{D}_PG &=\mathbf{D}_PA^*\mathbf{D}_PP+\mathbf{D}_PA\mathbf{D}_P-{\mathbf{D}_P}^2S\\
&=(S^*-P^*S)P+(S-S^*P)-(I-P^*P)S=0.\end{align*} Now $\langle
Gh,\mathbf{D}_Ph' \rangle=\langle \mathbf{D}_PGh,h' \rangle=0 \quad
\mbox{for all }h,h'\in\mathcal{H}$. This shows that $G=0$ and hence
$A^*\mathbf{D}_PP+A\mathbf{D}_P=\mathbf{D}_PS$. Therefore $T_AV_0=V_0T_A$. \\
Now
\begin{align*}
T_A^*V_0(h_0,h_1,h_2,\dots)&=T_A^*(Ph_0,\mathbf{D}_Ph_0,h_1,h_2,\dots)\\
&=(S^*Ph_0+\mathbf{D}_PA\mathbf{D}_Ph_0,A^*\mathbf{D}_Ph_0+Ah_1,A^*h_1+Ah_2,A^*h_2+Ah_3,\dots).
\end{align*}
Since $A$ is a solution of (\ref{e0}), we have
$S^*P+\mathbf{D}_{p}A\mathbf{D}_{p}=S$. Therefore we have $T_A^*V_0
=T_A$.

We now show that $r(T_A)\leq2$ which completes the proof. The
numerical radius of $A$ is not greater than $1$ by Theorem
\ref{greatlem}.\\
It is clear from the definition that $T_A$ has the matrix form
$$T_A= \begin{pmatrix} S&0&0&0&\dots\\
A^*\mathbf{D}_{p}&A&0&0&\dots\\
0&A^*&A&0&\dots\\
0&0&A^*&A&\dots\\
\dots&\dots&\dots&\dots&\dots\\
\end{pmatrix},$$ with respect to the decomposition $\mathcal{H}\oplus \mathcal{D}_p \oplus \mathcal{D}_p \oplus
\mathcal{D}_p\oplus...$ of $\mathcal{K}_0$. Again since $T_A=
\begin{pmatrix} S&0\\C&D
\end{pmatrix}$ on $\mathcal{H}\oplus
l^2(\mathcal{D}_P)= \mathcal{K}_0$, where $C=
\begin{pmatrix}
A^* \mathbf D_P\\0\\0\\ \vdots
\end{pmatrix} \mbox{ and }D=
\begin{pmatrix} A&0&0&\dots\\A^*&A&0&\dots\\0&A^*&A&\dots\\ \dots&\dots&\dots&\dots
\end{pmatrix}$, we have by Lemma 1 of \cite{hong} that $\sigma(T_A)\subseteq\sigma(S)\cup \sigma(D)$.
We shall be done if we show that $r(S)$ and $r(D)$ are not greater
than $2$. We show that $\|D\|\leq2.$ Let us define \begin{align*}
\varphi :& \,\mathbb{D}\rightarrow \mathcal{B}(\mathcal{D}_P)\\ & z
\rightarrow A+A^*z.
\end{align*}
Clearly $\varphi$ is holomorphic, bounded and continuous on the
boundary $\partial \mathbb{D}=\mathbb{T}$ of the disc. For $z=e^{-2i
\theta}\in\mathbb{T}$ we have \begin{align*} \|A+A^*z\|&= \|A+e^{-2i
\theta}A^*\|\\&=\|e^{i \theta}A+e^{-i \theta}A^*\|\\&=
\sup_{\|x\|\leq1}|\langle (e^{i \theta}A+e^{-i \theta}A^*)x,x
\rangle| \quad [\mbox{since} \;\; e^{i \theta}A+e^{-i \theta}A^*
\;\mbox{is self adjoint}]\\& \leq \omega(A)+\omega(A^*)\\& \leq 2.
\end{align*} Therefore by \textit{Maximum Modulus Principle}, $\|A+A^*z\|\leq2$
for all $z\in \overline{\mathbb{D}}$ and $\|\varphi\|\leq2$. Let
$$\widehat{A}_n={\begin{pmatrix}
A&0&0&\dots&0\\A^*&A&0&\dots&0\\0&A^*&A&\dots&0\\
\dots&\dots&\dots&\dots&\dots\\0&0&\dots&A^*&A
\end{pmatrix}}_{n\times n}\quad \mbox{on}\;\; \displaystyle \underbrace{\mathcal{D}_P
\oplus \mathcal{D}_P\oplus \dots \oplus \mathcal{D}_P}_{n-times}=
\mathbf{E}_n.$$ Let $f= \displaystyle \bigoplus_0^{n-1}f_i$ and
$g=\displaystyle\bigoplus_0^{n-1}g_i$ be two arbitrary elements in
$\mathbf{E}_n$. Let us consider the polynomials
$p(z)=\displaystyle\sum_{i=0}^{n-1}z^if_i$ and
$q(z)=\displaystyle\sum_{i=0}^{n-1}z^ig_i$ with values in
$\mathcal{D}_P$. Now
\begin{align*} |\langle \widehat{A}_nf,g
\rangle_{\mathbf{E}_n}|= |\langle
\widehat{A}_n(\displaystyle\bigoplus_0^{n-1}f_i),(\displaystyle
\bigoplus _0^{n-1}g_i)
\rangle_{\mathbf{E}_n}|&=|\frac{1}{2\pi}\int_0^{2\pi}\langle \phi
(e^{it})p(e^{it}),q(e^{it}) \rangle_{\mathcal{D}_P} \,dt|\\&
\leq\|\phi(e^{it})p(e^{it}) \|_{L^2} \|q(e^{it})\|_{L^2}\\ & \leq
2\|p(e^{it}) \|_{L^2} \| q(e^{it}) \|_{L^2} \quad [\mbox{since }
\|\phi\|\leq2]\\& =2\|\displaystyle
\bigoplus_0^{n-1}f_i\|_{\mathbf{E}_n} \|\displaystyle
\bigoplus_0^{n-1}g_i\|_{\mathbf{E}_n}\\&=2\|f\|\|g\|
\end{align*}
This implies that $\|\widehat{A}_n \|\leq2.$ Now we define $D_n$ on
$\mathbf{E}_n \oplus \mathbf{E}_{\infty}=l^2(\mathcal{D}_P)$, where
$\mathbf{E}_{\infty}=l^2(\mathcal{D}_P)\ominus \mathbf{E}_n$, as
$D_n=\begin{pmatrix} \widehat{A}_n&0\\0&0
\end{pmatrix}$. Then $\|D_n\|=\|\widehat{A}_n\|\leq2$ and $D_n \rightarrow
D$ strongly as $n\rightarrow \infty.$ Hence $\|D\|\leq2.$ Again
since $(S,P)$ is a $\Gamma$-contraction, $r(S)\leq \|S\|\leq2$.
Since both of $r(S),r(D)$ are not greater than $2$, $r(T_A)\leq2$. Hence done.\\
\textbf{(2)} Obviously $V_0=\begin{pmatrix} P&0\\C_1&D_1
\end{pmatrix}$ with respect to the decomposition $\mathcal H\oplus l^2(\mathcal {D}_P)$
of $\mathcal K_0$, where $$C_1=\begin{pmatrix} \mathbf{D}_P
\\0\\0\\ \vdots \end{pmatrix} \mbox{ from } \mathcal H \rightarrow
\mathcal{D}_P\oplus\mathcal{D}_P\oplus\mathcal{D}_P\oplus\dots
\mbox{ and } D_1=\begin{pmatrix}
0&0&0&\dots\\I&0&0&\dots\\0&I&0&\dots\\\dots&\dots&\dots&\dots
\end{pmatrix} \mbox{ on }
\mathcal{D}_P\oplus\mathcal{D}_P\oplus\mathcal{D}_P\oplus\dots.$$
Since $(\widehat{T},V_0)$ on $\mathcal K_0$ is a $\Gamma$-isometric
dilation of $(S,P)$, we have $\widehat{T}^*|_{\mathcal H}=S^*$ and
$V_0^*|_{\mathcal H}=P^*$. Therefore $\widehat{T}$ on $\mathcal
H\oplus l^2(\mathcal D_P)$ has matrix form $\widehat{T}=
\begin{pmatrix} S&0\\E&F
\end{pmatrix}$. Let us define \begin{align*}U_1: & H^2(\mathcal D_P)\rightarrow \mathcal D_P\oplus \mathcal D_P\oplus \mathcal
D_P\oplus\dots\\& z^n\mapsto (\displaystyle
\underbrace{0,0,\dots,0}_n,1,0,0,\dots).
\end{align*} The action of $U_1$ on an arbitrary vector is clear from its action on the basis
$\{1,z,z^2,\dots\}$ of $H^2(\mathcal D_P)$. Since it maps a basis of
$H^2(\mathcal D_P)$ to a basis of $\mathcal D_P\oplus \mathcal
D_P\oplus \mathcal D_P\oplus\dots$ in a one-to-one fashion, $U_1$ is
a unitary operator. Therefore the spaces $\mathcal K_0$ and
$\mathcal H \oplus H^2(\mathcal D_P)$ are isomorphic. Let $U=U_1^*$.
Then $\widehat{T}$ and $V_0$ on $\mathcal K_0$ are respectively
identified with the operators
$$\widetilde{T}=\begin{pmatrix} S&0\\UE&UFU^* \end{pmatrix} \mbox{ and } \widetilde{V_0}=\begin{pmatrix} P&0\\UC_1&UD_1U^* \end{pmatrix}
\mbox{ on } \mathcal H\oplus H^2(\mathcal D_P).$$ Therefore
$(\widetilde{T},\widetilde{V_0})$ is a $\Gamma$-isometric dilation
of $(S,P)$. We now show that $UD_1U^*$ is same as the multiplication
operator $M_z^{\mathcal D_P}$ on $H^2(\mathcal D_P)$. For a basis
vector $z^n$ of $H^2(\mathcal D_P)$ we have
\begin{align*}UD_1U^*(z^n)&=U \begin{pmatrix}
0&0&0&0&\dots\\I&0&0&0&\dots\\0&I&0&0&\dots\\0&0&I&0&\dots
\\\dots&\dots&\dots&\dots&\dots \end{pmatrix} \begin{pmatrix}
0\\\vdots\\0\\1\\0\\\vdots \end{pmatrix}, \quad 1 \mbox{ at }
(n+1)\mbox{th place}\\ &=U
\begin{pmatrix} 0\\\vdots\\0\\1\\0\\\vdots \end{pmatrix}, \quad 1\mbox{ at }(n+2)\mbox{th place}\\& =z^{n+1}=M_z^{\mathcal D_P}(z^n).\end{align*}
Hence $UD_1U^*=M_z^{\mathcal D_P}$. By the commutativity of
$\widetilde{T}$ and $\widetilde{V}_0$ we have the commutativity of
$UFU^*$ and $UD_1U^*(=M_z^{\mathcal D_P})$. Therefore
$UFU^*=M_{\varphi}^{\mathcal D_P}$ for some $\varphi\in
H^{\infty}(\mathcal B (\mathcal D_P))$. Thus
$$\widetilde{T}=\begin{pmatrix}S&0\\UE&M_{\varphi}^{\mathcal D_P} \end{pmatrix}\mbox{ and }
\widetilde{V_0}=\begin{pmatrix} P&0\\UC_1&M_z^{\mathcal D_P}
\end{pmatrix} \mbox{ on } \mathcal H\oplus H^2(\mathcal D_P).$$By
$\widetilde{T}=\widetilde{T}^*\widetilde{V_0}$, we get
$$
\begin{pmatrix} S&0\\UE&M_{\varphi}^{\mathcal D_P} \end{pmatrix}=
\begin{pmatrix} S^*&E^*U^*\\0&{M_{\varphi}^{\mathcal D_P}}^* \end{pmatrix}
\begin{pmatrix} P&0\\UC_1&M_z^{\mathcal D_P} \end{pmatrix} =
\begin{pmatrix} S^*P+E^*C_1&E^*U^*M_z^{\mathcal D_P}\\{M_{\varphi}^{\mathcal D_P}}^*UC_1&{M_{\varphi}^{\mathcal D_P}}^*M_z^{\mathcal D_P}
\end{pmatrix},
$$ which gives \begin{eqnarray}\begin{cases} \label{e9}
&(\mbox{i})\; S-S^*P=E^*C_1\\&(\mbox{ii})\;
UE={M_{\varphi}^{\mathcal D_P}}^*UC_1\\&(\mbox{iii})\;
M_{\varphi}^{\mathcal D_P}={M_{\varphi}^{\mathcal D_P}}^*M_z.
\end{cases}\end{eqnarray}
 From (\ref{e9})-(iii), it is clear by considering the power
 series expansion that $\varphi(z)=A_0+A_0^*z$, for some $A_0\in\mathcal B(\mathcal
 D_P)$. We now show that if $D_0=\begin{pmatrix} A_0&0&0&\dots\\A_0^*&A_0&0&\dots\\0&A_0^*&A_0&\dots\\\dots&\dots&\dots&\dots \end{pmatrix}$
on $\mathcal D_P\oplus \mathcal D_P\oplus \mathcal D_P\oplus\dots$,
then $UD_0U^*=M_{\varphi}^{\mathcal D_P}$. For a basis vector $z^n$
of $H^2(\mathcal D_P)$ we have
$$
UD_0U^*(z^n)=U\begin{pmatrix}A_0&0&0&\dots\\A_0^*&A_0&0&\dots\\0&A_0^*&A_0&\dots\\
\dots&\dots&\dots&\dots
\end{pmatrix}
\begin{pmatrix}
0\\\vdots\\0\\1\\0\\\vdots
\end{pmatrix}
=U
\begin{pmatrix}
0\\\vdots\\0\\A_0(1)\\A_0^*(1)\\0\\\vdots
\end{pmatrix}
=A_0(1)z^n+A_0^*(1)z^{n+1}={M_{A_0+A_0^*z}^{\mathcal D_P}}(z^n)
$$
Thus $UD_0U^*=M_{\varphi}^{\mathcal D_P}$ and hence $F=D_0$.
Combining this with (\ref{e9})-(ii), we get
$UE={M_{\varphi}^{\mathcal D_P}}^*UC_1=UD_0^*U^*UC_1=UD_0^*C_1$,
i.e, $E=D_0^*C_1$. Therefore
$$\widehat{T}=\begin{pmatrix}
S&0\\D_0^*C_1&D_0
\end{pmatrix} \mbox{ on } \mathcal H\oplus l^2(\mathcal D_P).$$
Considering the above stated matrix forms of $D_0$ and $C_1$ we get
$D_0^*C_1=\begin{pmatrix} A_0^* \mathbf D_P
\\0\\0\\\vdots \end{pmatrix}.$ Hence with respect to the
decomposition $\mathcal H\oplus \mathcal D_P\oplus\mathcal
D_P\oplus\dots$ of $\mathcal K_0$, we have
$$\widehat{T}=\begin{pmatrix} S&0&0&0&\dots\\A_0^*\mathbf D_P &A_0&0&0&\dots\\
0&A_0^*&A_0&0&\dots\\0&0&A_0^*&A_0&\dots\\\dots&\dots&\dots&\dots&\dots
\end{pmatrix}.$$ Also by (\ref{e9})-(i),
\begin{align*}
S-S^*P=E^*C_1 &=C_1^*D_0C_1\\&=\begin{pmatrix} \mathbf D_P&0&0&\dots
\end{pmatrix}
\begin{pmatrix}
A_0&0&0&\dots\\A_0^*&A_0&0&\dots\\0&A_0^*&A_0&\dots\\
\dots&\dots&\dots&\dots
\end{pmatrix}
\begin{pmatrix}
\mathbf D_P\\0\\0\\\vdots
\end{pmatrix}=\mathbf D_PA_0 \mathbf D_P,
\end{align*} which shows that $A_0$ satisfies the fundamental
equation (\ref{e0}). By uniqueness of solution, $A=A_0$ and hence
$\widehat{T}=T_A$.

\textbf{(3)} Let $(T,V)$ defined on $\mathcal K$ be a minimal
isometric dilation of $(S,P)$, where $V$ is a minimal isometric
dilation of $P$. Since $V$ on $\mathcal K$ is a minimal isometric
dilation of $P$, there is a unitary $$U: \mathcal K\rightarrow
\mathcal K_0 (=\mathcal H\oplus \mathcal{D}_P\oplus
\mathcal{D}_P\oplus \dots)$$ such that $UVU^*=V_0$. Let $T^{\flat}
=UTU^*$. Then $(T^{\flat},V_0)$ on $\mathcal K_0$ is a
$\Gamma$-isometry dilation of $(S,P)$. Therefore by part-($2$),
$T^{\flat} =T_A$ and consequently $(T,V)$ is unitarily equivalent to
$(T_A,V_0)$.
\end{proof}

As a consequence of the dilation theorem above, we have a new and
elegant characterization for $\Gamma$-contractions.

\begin{thm}
Let $(S,P)$ be a commuting pair of operators defined on $\mathcal
H$. Then $(S,P)$ is a $\Gamma$-contraction if and only if spectral
radius of $S$ is not greater than $2$ and the fundamental equation
$S-S^*P=\mathbf{D}_PX\mathbf{D}_P$ has a solution $A$ with
$\omega(A)\leq1$.
\end{thm}
\begin{proof}
Let there be a solution $A$ to the fundamental equation
$S-S^*P=\mathbf{D}_PX\mathbf{D}_P$ with $\omega(A)\leq1$ for such a
pair $(S,P)$. Then by the dilation theorem (Theorem
\ref{schaeffer}), we can construct a $\Gamma$-isometry $(T_A,V_0)$
of $(S,P)$. Now clearly $(S,P)$ can be recovered by compressing
$(T_A,V_0)$ to the common co-invariant subspace $\mathcal H$. So
$(S,P)$ is a $\Gamma$-contraction.

The converse is just the Theorem \ref{greatlem}.
\end{proof}
\begin{rem}
{\em We call the unique solution $A$ of the operator equation
(\ref{e0}) for a $\Gamma$-contraction $(S,P)$, the fundamental
operator of $(S,P)$}.
\end{rem}

We now give another explicit construction of a $\Gamma$-isometric
dilation of a pure $\Gamma$-contraction. This is very convenient to
reap some beautiful consequences.

\begin{thm} \label{main}
Let $(S,P)$ be a $\Gamma$-contraction on a Hilbert space
$\mathcal{H}$ where $P$ is in $C_{\cdot0}$. Let $B$ be the solution
of the fundamental equation (\ref{e1}). Let us consider the
operators $T,V$ on $\mathcal{E}=H^2(\mathbb{D}) \otimes
\mathcal{D}_{P^*}$ defined as
$$T=I\otimes B^*+M_z\otimes B \mbox{ and } V=M_z\otimes I.$$ Then
$(T,V)$ is a $\Gamma$-isometric dilation of $(S,P)$.
\end{thm}

\begin{proof}
Since $B$ is the solution of the equation (\ref{e1}), by Theorem
\ref{greatlem}, the numerical radius of $B$ is not greater than one.
In order to prove that $(T,V)$ is a $\Gamma$-isometric dilation of
$(S,P)$ we shall show the following steps:
\begin{enumerate}
\item the pair (T,V) is a $\Gamma$-{\em isometry} on
    $\mathcal{E}$.
    \item The space $\mathcal{H}$ can be thought of as a
        subspace of $\mathcal{E}$, i.e, there is an isometric
        embedding of $\mathcal H$ in $\mathcal{E}$.
\item After identification of $\mathcal H$ with this isometric
    image, $V^*\mathcal H \subseteq \mathcal{H}$ and
    $V^*|_{\mathcal{H}}=P^*$. Also, $T^*\mathcal{H}\subseteq
    \mathcal{H} $ and $T^*|_{\mathcal{H}}=S^*.$
\end{enumerate}

$V$ is clearly an isometry (it is a shift of some multiplicity) and
obviously it commutes with $T$. Also
$$T =(I\otimes B^*)+(I\otimes B)(M_z\otimes I)= C+C^*V,$$ where $C=I\otimes
B^*$. Obviously $C$ and $C^*$ commute with $V$ and $\omega(C)\leq1$.
Therefore by Theorem \ref{goodthm}, $(T,V)$ is a $\Gamma$-isometry.

Now we embed $\mathcal H$ isometrically inside $H^2 \otimes
{\mathcal D}_{P^*}$ by defining $W:\mathcal{H} \rightarrow
\mathcal{E}$ as
 $ h\mapsto \sum_{n=0}^{\infty} z^n\otimes \mathbf{D}_{P^*}{P^*}^n h$.  \begin{align*}
\|Wh\|^2 &= \|\displaystyle \sum_{n=0}^{\infty}{z^n\otimes
\mathbf{D}_{P^*}{P^*}^n h}\|^2 \\&= \langle \displaystyle
\sum_{n=0}^{\infty}{z^n\otimes \mathbf{D}_{P^*}{P^*}^n
h}\;,\;\displaystyle \sum_{m=0}^{\infty}{z^m\otimes
\mathbf{D}_{P^*}{P^*}^m h} \rangle \\& = \displaystyle
\sum_{m,n=0}^{\infty} \langle z^n,z^m \rangle \langle
\mathbf{D}_{P^*}{P^*}^nh\;,\;\mathbf{D}_{P^*}{P^*}^mh \rangle \\&
=\displaystyle \sum_{n=1}^{\infty}{\langle P^n\mathbf{D}_{P^*}^2
{P^*}^nh,h \rangle}\\&= \displaystyle \sum_{n=0}^{\infty}\langle
P^n(I-PP^*){P^*}^nh,h \rangle\\& = \displaystyle
\sum_{n=0}^{\infty}\{\langle P^n{P^*}^nh,h \rangle-\langle
P^{n+1}{P^*}^{n+1}h,h \rangle\} \\&= \|h\|^2-\lim_{n \rightarrow
\infty}\|{P^*}^nh\|^2.
\end{align*} Since $P\in C_{\cdot 0},\; \displaystyle \lim_{n\rightarrow
\infty}\|{P^*}^nh\|^2=0$ and hence $\|Wh\|=\|h\|.$ Therefore $W$ is
an isometry. Let $L=W^*.$

For a basis vector $z^n\otimes \xi$ of $\mathcal{E}$ we have
$$ \langle L(z^n\otimes \xi),h \rangle  = \langle z^n \otimes \xi ,
\displaystyle \sum_{k=0}^{\infty}{z^k \otimes
\mathbf{D}_{P^*}{P^*}^kh} \rangle  = \langle \xi ,
\mathbf{D}_{P^*}{P^*}^nh \rangle  = \langle P^n \mathbf{D}_{P^*}\xi
, h \rangle. $$ This implies that $$L(z^n \otimes \xi)=P^n
\mathbf{D}_{P^*} \xi, \quad for\; n=0,1,2,3,...$$ Therefore
$$ \langle L(M_z \otimes I)(z^n \otimes \xi), h \rangle = \langle
z^{n+1} \otimes \xi , \displaystyle \sum_{k=0}^{\infty}{z^k \otimes
\mathbf{D}_{P^*}{P^*}^kh} \rangle  = \langle \xi,
\mathbf{D}_{P^*}{P^*}^{n+1}h \rangle = \langle
P^{n+1}\mathbf{D}_{P^*} \xi,h \rangle. $$ Consequently, $LV = PL$ on
vectors of the form $z^n \otimes \xi$ which span $H^2 \otimes
{\mathcal D}_{P^*}$ and hence
\begin{eqnarray}\label{e7} LV = PL.\end{eqnarray} Therefore $V^*$ leaves the range of $L^*$
(isometric copy of $\mathcal H$) invariant and
$V^*|_{L^*\mathcal{H}}=L^*P^*L$ which is the isometric copy of the
operator $P^*$ on range of $L^*$. For the next step,
\begin{align*} LT(z^n \otimes \xi)=L(I \otimes B^*+M_z \otimes
B)(z^n \otimes \xi) &=L(I \otimes B^*)(z^n \otimes \xi)+L(M_z
\otimes B)(z^n \otimes \xi)\\&=L(z^n \otimes B^*\xi)+L(z^{n+1}
\otimes B \xi)\\&=P^n\mathbf{D}_{P^*}B^*\xi +
P^{n+1}\mathbf{D}_{P^*}B\xi.
\end{align*} Again $SL(z^n \otimes \xi)=SP^n \mathbf{D}_{P^*} \xi$.
Therefore for showing $LT=SL$ it is enough to show that
\begin{align*}
& P^n\mathbf{D}_{P^*}B^*+P^{n+1}\mathbf{D}_{P^*}B=SP^n\mathbf{D}_{P^*}=P^nS\mathbf{D}_{P^*}\\
& \mbox{i.e,}\;
\mathbf{D}_{P^*}B^*+P\mathbf{D}_{P^*}B=S\mathbf{D}_{P^*}.
\end{align*}
Let $H=\mathbf{D}_{P^*}B^*+P\mathbf{D}_{P^*}B-S\mathbf{D}_{P^*}$.
Then $H=0$ by an argument similar to the one given in the proof of
Theorem \ref{schaeffer} to show that $G=0$. So we have
$$\mathbf{D}_{P^*}B^*+P\mathbf{D}_{P^*}B=S\mathbf{D}_{P^*}$$ and hence
\begin{eqnarray} \label{e8}L(I \otimes B^*+M_z \otimes
B)=SL\end{eqnarray} which is similar to the equation (\ref{e7}).
This shows that $T^*$ leaves $L^*(\mathcal H)$ invariant as well as
$T^*|_{L^*(\mathcal H)}=L^*S^*L$. Hence we are done.
\end{proof} \vspace{0.10in}

\begin{rem}
{\em In particular when $\|P\|<1$ the unique solutions $A$ of
(\ref{e0}) and $B$ of (\ref{e1}) coincide with
$\mathbf{D}_P^{-1}(S-S^*P)\mathbf{D}_P^{-1}$ and
$\mathbf{D}_{P^*}^{-1}(S^*-SP^*)\mathbf{D}_{P^*}^{-1}$
respectively.}
\end{rem}
\begin{cor} If $(S,P)$ is a $\Gamma$-contraction with $P\in C_{\cdot 0}$,
then $S=C+PC^*$ for some $C$ with $\omega(C)\leq1.$ \end{cor}
\begin{proof}
By the previous theorem, if $(T,V)$ is a $\Gamma$-isometric dilation
of $(S,P)$ from (\ref{e8}) we have \begin{align*} & LT=L(I \otimes
B^*+M_z \otimes B)=SL \\&\mbox{ or } L(I \otimes B^*+M_z \otimes
B)L^*=S,\; \mbox{since }L^* \mbox{ is isometry }
\\& \mbox{ or } L(I \otimes B^*)L^*+L(M_z \otimes B)L^*=S\\&\mbox{
or } L(I \otimes B^*)L^*+L(M_z \otimes I)(I \otimes B)L^*=S \\&
\mbox{ or } L(I \otimes B^*)L^*+PL(I \otimes B)L^*=S, \; \mbox{
since }L(M_z \otimes I)=PL.
\end{align*} Taking $C= L(I \otimes B^*)L^*$ we get the stated
form of $S,$ and $\omega(C)\leq1$ is obvious.
\end{proof}\vspace{0.07in}

\begin{obs}
If $(S,P)$ is a $\Gamma$-contraction with $P\in C_{\cdot 0}$, then
$S$ can also have the form $S=C_1+C_1^*P,$ where $\omega(C_1)\leq1.$
\end{obs}
\begin{proof}
Clearly $(S^*,P^*)$ is also a $\Gamma$-contraction and by the
previous result, $S^*=C+P^*C^*$ for some $C$ with $\omega(C)\leq1.$
This implies that $S=C^*+CP=C_1+C_1^*P$ where $C_1=C^*$.
\end{proof}
\begin{obs}
If $(S,P)$ is a $\Gamma$-contraction with $\|P\|<1$, then there is a
unique $C$ such that $S=C+C^*P$.
\end{obs}
\begin{proof}
Let there be $C_1$ and $C_2$ such that $S=C_1+C_1^*P$ and
$S=C_2+C_2^*P$. Then we have $C+C^*P=0$, where $C=C_1-C_2$. Now
$$\|C\|=\|-C^*P\|\leq \|C\|\|P\|< \|C\| \quad \mbox{as }
\|P\|<1.$$ This shows that $C=0$ and consequently $C_1=C_2$.
\end{proof}

For a polynomially convex compact subset $X$ of $\mathbb C^d$ and a
tuple of commuting bounded operators $\underline{A}=(A_1,\dots,A_d)$
on a Hilbert space $\mathcal H$, a normal $\partial X$-dilation
$\underline{N}=(N_1,\dots,N_d)$ is a tuple of commuting bounded
operators on a Hilbert space $\mathcal K\supseteq \mathcal H$ such
that the Taylor joint spectrum $\sigma_T(\underline{N})\subseteq
\partial X$ and $$ p(\underline{A})=P_{\mathcal
H}p(\underline{N})|_{\mathcal H}, \quad \textup{ for any } p\in
\mathbb C[z_1,\dots,z_d].
$$
It is clear that if $\underline{A}$ has a normal $\partial
X$-dilation, then $X$ is a spectral set for $\underline{A}$. In
general, it is difficult to determine the converse, i.e, if $X$ is a
spectral set for $A$ then whether or not $\underline{A}$ has a
normal $\partial X$-dilation. It was shown by Agler and Young that a
pair of commuting bounded operators $(S,P)$ has $\Gamma$ as a
spectral set if and only if it has a normal $\partial X$-dilation.
One of the contributions of this paper has been to add that $\Gamma$
is a spectral set for a commuting pair $(S,P)$ if and only if the
fundamental equation for $(S,P)$ can be solved with a solution of
numerical radius not greater than one.

\textbf{Acknowledgement.} We are thankful to Professor Gadadhar
Misra for stimulating conversations.

\end{document}